\documentclass[aos,MSNbibl,dvips]{arximspdf}

\doi{10.1214/13-AOS1133} 
\volume{41}
\issue{4}
\pubyear{2013}
\firstpage{1999}
\lastpage{2028}

\makeatletter
\newcommand{\rrvert}{\vert}
\newcommand{\llvert}{\vert}

\newtheorem{theorem}{Theorem}

\newproclaim{condition}{Condition}
\newproclaim{definition}{Definition}

\newtheorem{corollary}{Corollary}

\newproclaim{remark}{Remark}

\newcommand{\al}{\alpha}
\newcommand{\be}{\beta}
\newcommand{\ga}{\gamma}
\newcommand{\la}{\lambda}
\newcommand{\te}{\theta}
\newcommand{\ta}{\tau}
\newcommand{\veps}{\varepsilon}
\newcommand{\vphi}{\varphi}
\newcommand{\pli}{+\infty}

\newcommand{\cJ}{\mathcal{J}}
\newcommand{\cN}{\mathcal{N}}

\newcommand{\RR}{\mathbb{R}}

\newcommand{\given}{|}

\newcommand{\bfi}{\bar{\Phi}}

\newcommand{\tlk}{\theta_{lk}}
\newcommand{\tolk}{\theta_{0,lk}}
\newcommand{\sil}{\sigma_{l}}
\newcommand{\elk}{\veps_{lk}}

\newcommand{\rn}{\sqrt{n}}

\newcommand{\ix}{\mathbb{X}}
\newcommand{\dob}{\mathbb{W}}

\makeatother

\begin{document}
\begin{frontmatter}

\title{Nonparametric Bernstein--von Mises theorems in Gaussian white noise}
\runtitle{Nonparametric BvMs}

\begin{aug}
\author[A]{\fnms{Isma\"el} \snm{Castillo}\thanksref{t1}\ead[label=e1]{ismael.castillo@upmc.fr}}
\and
\author[B]{\fnms{Richard} \snm{Nickl}\corref{}\ead[label=e2]{r.nickl@statslab.cam.ac.uk}}
\runauthor{I. Castillo and R. Nickl}
\affiliation{CNRS and University of Cambridge}
\address[A]{CNRS\\
Laboratoire Probabilit\'es \\
\quad et Mod\`eles Al\'eatoires\\
LPMA\\
Universit\'es Paris VI \& VII\\
B\^atiment Sophie Germain \\
75205 Paris Cedex 13\\
France\\
\printead{e1}} 
\address[B]{Statistical Laboratory\\
Department of Pure Mathematics \\
\quad and Mathematical Statistics\\
University of Cambridge\\
CB3 0WB Cambridge\\
United Kingdom \\
\printead{e2}}
\end{aug}

\thankstext{t1}{Supported in part by ANR Grants ``Banhdits''
ANR-2010-BLAN-0113-03 and ``Calibration'' ANR-2011-BS01-010-01.}

\received{\smonth{1} \syear{2013}}
\revised{\smonth{5} \syear{2013}}

%
\begin{abstract}
Bernstein--von Mises theorems for nonparametric Bayes priors in the
Gaussian white noise model are proved. It is demonstrated how such
results justify Bayes methods as efficient frequentist inference
procedures in a variety of concrete nonparametric problems.
Particularly Bayesian credible sets are constructed that have
asymptotically exact $1-\alpha$ frequentist coverage level and whose
$L^2$-diameter shrinks at the minimax rate of convergence (within
logarithmic factors) over H\"older balls. Other applications include
general classes of linear and nonlinear functionals and credible bands
for auto-convolutions. The assumptions cover nonconjugate product
priors defined on general orthonormal bases of $L^2$ satisfying weak
conditions.
\end{abstract}

%
\begin{keyword}[class=AMS]
\kwd[Primary ]{62G20}
\kwd[; secondary ]{62G15}
\kwd{62G08}
\end{keyword}
\begin{keyword}
\kwd{Bayesian inference}
\kwd{plug-in property}
\kwd{efficiency}
\end{keyword}

\end{frontmatter}

\section{Introduction}

Consider observing a random sample $X^{(n)}$ of size $n$, or at noise
level $n^{-1/2}$, drawn from distribution $P_{f}^n$, and indexed by
some unknown parameter $f \in\mathcal F$. The Bayesian paradigm views
the sample as having law $P_{f}^n$ conditionally on $f$, that is,
$X^{(n)}|f \sim P_{f}^n$, and the law of $f$ is a \textit{prior
probability distribution} $\Pi$ on some $\sigma$-field $\mathcal B$
of $\mathcal F$. The random variable $f|X^{(n)}$ then has a law on
$\mathcal F$ which is known as the \textit{posterior distribution},
denoted by $\Pi(\cdot|X^{(n)})$.
Bayesian inference on $f$ is then entirely based on this posterior
distribution---it gives access to point estimates for $f$, credible
sets and tests in a natural way.

It is of interest to analyse the behaviour of $\Pi(\cdot|X^{(n)})$
under the frequentist sampling assumption that $X^{(n)}$ is drawn from
$P_{f_0}^n$ for some fixed nonrandom $f_0 \in\mathcal F$. If $\mathcal
F$ is a \textit{finite-dimensional} space, then posterior-based
inference has a fundamental justification through the Bernstein--von
Mises (BvM) theorem, first discovered by Laplace \cite{L1810},
developed by von Mises \cite{vM31}, and put into the framework of
modern parametric statistics by Le Cam \cite{LC86}. It states that,
under mild and universal assumptions on the prior, the posterior
distribution approximately equals a normal distribution $N(\hat f_n,
i(f))$ centered at an efficient estimator $\hat f_n$ for $f$ and with a
covariance $i(f)$ that attains the Cram\'er--Rao bound in the
statistical model considered: as $n \to\infty$
%
\begin{equation}
\label{bvm0} \sup_{B \in\mathcal B} \bigl\llvert\Pi\bigl(B|X^{(n)}
\bigr) - N\bigl(\hat f_n, i(f_0)\bigr) (B) \bigr\rrvert
\to^{P_{f_0}^n} 0.
\end{equation}
As a consequence posterior-based inference asymptotically coincides
with inference based on standard efficient, $1/\sqrt n$-consistent
frequentist estimators $\hat f_n$ of $f$, giving a rigorous asymptotic
justification of Bayesian methods.

The last decade has seen remarkable activity in the development of
\textit{nonparametric} Bayes procedures, where $\mathcal F$ is taken to
be an infinite-dimensional space, typically consisting of functions or
infinite vectors: nonparametric regression, classification, density
estimation, normal means, and Gaussian white noise models come to mind,
and a variety of nonparametric priors have been devised in the
literature for such models. Posteriors can be computed efficiently by
algorithms such as MCMC, and they provide broadly applicable Bayesian
inferential tools for nonparametric problems. It is natural to ask
whether an analogue of (\ref{bvm0}) can still be proved in such
situations, as it would give a general justification for the use of
nonparametric Bayes procedures. Although remarkable progress has been
made in the understanding of the frequentist properties of
nonparametric Bayes procedures (we refer here only to some of the key
papers such as \cite{GGV00,SW01,GV07,VV08,VV09} and references
therein), a fully satisfactory answer to the BvM-question seems not to
have been found. A first reason is perhaps that it is not immediately
clear what $N(\hat f_n, i(f_0))$ should be replaced by in the
infinite-dimensional situation---Gaussian distributions over
infinite-dimensional spaces $\mathcal F$ are much more complex objects,
and their existence in the form relevant here depends on the topology
that $\mathcal F$ is endowed with. Another reason is that the commonly
used loss functions in nonparametric statistics (such as $L^p$-type
loss) do not admit $1/\sqrt n$-consistent estimators---the LAN-type
local approximations of the likelihood function used in the proof of
the finite-dimensional BvM theorem are thus not accurate enough in such
metrics.

One way around these problems is to weaken the loss function on
$\mathcal F$ so that $1/\sqrt n$-consistent estimation with Gaussian
limits is possible even in nonparametric models. For example, in the
situation where one observes $X_1, \ldots, X_n \sim^{\mathrm{i.i.d.}} P$ on
$[0,1]$, and $P_n = n^{-1} \sum_{i=1}^n \delta_{X_i}$ is the
empirical measure estimating $P$, then for any $P$-Donsker class
$\mathcal H$ of functions $h\dvtx  [0,1] \to\mathbb R$,
%
\begin{equation}
\label{pipu} \sup_{h \in\mathcal H}\biggl\llvert\int
hd(P_n-P) \biggr\rrvert= \sup_{h \in
\mathcal H}\Biggl\llvert
\frac{1}{n} \sum_{i=1}^n
\bigl(h(X_i)-Eh(X_1)\bigr) \Biggr\rrvert=
O_P\bigl(n^{-1/2}\bigr),
\end{equation}
in fact $P_n$ is an efficient estimator for $P$ in the space $l^\infty
(\mathcal H)$ of bounded functions on $\mathcal H$, attaining the
Cram\'er--Rao information bound for the fully nonparametric model; cf.
3.1.11 in \cite{VW96}. Hence one may try to prove a BvM-type result by
endowing the parameter space $\mathcal F$ with the loss function coming
from an $l^\infty(\mathcal H)$-type space. The purpose of the present
paper is to investigate this approach rigorously in the setting of the
Gaussian white noise model, and with $\mathcal H$ a ball in a suitable
Sobolev space defined below. This makes the mathematical analysis
tractable without any severe loss of conceptual generality; see below
for a discussion of extensions to other models. Our main results will
imply that for a large and relevant class of nonparametric product
priors $\Pi$ that satisfy mild assumptions, and which do not require
conjugacy, one has
%
\begin{equation}
\label{bvm1} \sup_{A \in\mathcal A} \bigl\llvert\Pi\bigl(A|X^{(n)}
\bigr) - \mathcal N\bigl(\hat f_n, i(f_0)\bigr) (A) \bigr
\rrvert\to^{P_{f_0}^n} 0,
\end{equation}
where $\mathcal N$ is a Gaussian measure on $l^\infty(\mathcal H)$
centered at an efficient estimator $\hat f_n$ of~$f$, both to be
defined in a precise manner, and where the classes $\mathcal A$ consist
of measurable subsets of $l^\infty(\mathcal H)$ that have uniformly
smooth boundaries for the measure $\mathcal N$. We note that some
restrictions on the class $\mathcal A$ are necessary as one can show
that in the infinite-dimensional situation the Bernstein--von Mises
theorem cannot hold uniformly in all Borel sets of $l^\infty(\mathcal
H)$; see after Definition \ref{bvmp} for discussion. Our assumptions
apply to priors that produce posteriors which achieve frequentist
optimal contraction rates in stronger loss functions (such as
$L^2$-distance) and which resemble the state of the art prior choices
in the nonparametric Bayes literature.

Our abstract results only gain relevance through the fact that we can
demonstrate their applicability: the general result (\ref{bvm1}) will
be shown to imply that posterior-based credible regions give
asymptotically exact frequentist confidence sets in a variety of
concrete problems of nonparametric inference. A first important
application is to weighted $L^2$-ellipsoid credible regions for the
unknown parameter $f$, which are shown to have optimal width
$O_P(n^{-1/2})$ in $\ell^\infty(\mathcal H)$-loss and which
simultaneously are confidence sets that shrink in $L^2$-diameter at the
minimax rate (within log-factors) over H\"older balls. We further give
semiparametric applications to estimation of linear and nonlinear
functionals defined on $\mathcal F$, and to credible bands for
estimating an auto-convolution $f \ast f$.

A key point in these applications is related to the notion of the
``plug-in property'' coined by Bickel and Ritov \cite{BR03}. A
nonparametric estimator that is rate-optimal in a standard loss
function (such as $L^p$-loss) is said to have the plug-in property if
it simultaneously efficiently estimates, at $1/\sqrt n$-rate, a large
class of linear functionals. Standard frequentist estimators such as
kernel, wavelet and nonparametric maximum likelihood estimators satisfy
this property; in fact, one can even prove a corresponding uniform
central limit theorem in $l^\infty(\mathcal H)$ for such estimators;
see Kiefer and Wolfowitz \cite{KW76}, Nickl \cite{N07}, Gin\'{e} and
Nickl \cite{GN08,GN09a}. Our results imply that this is also true in
the Bayesian situation: the posterior contracts at the optimal rate in
$L^2$-loss and at the same time satisfies a Bernstein--von Mises
theorem in $l^\infty(\mathcal H)$.

Our techniques of proof rely on the structure of the Gaussian white
noise model and apply as well, with simple modifications, to fixed
design nonparametric regression settings. Our proofs moreover indicate
a strategy to obtain BvM results of this kind in general nonparametric
sampling models: the idea is to obtain ``semiparametric'' BvM-results for
many fixed linear functionals simultaneously, and to re-construct
``nonparametric'' norms from these functionals. At least for density
estimation such ideas can be shown to work as well, and this is subject
of forthcoming research.

There is important work on the BvM phenomenon for nonparametric
procedures that needs mentioning. Cox \cite{C93} and Freedman \cite
{Free99} have shown the impossibility of a nonparametric BvM result in
a strict $L^2$-setting. Leahu \cite{L11} derives interesting results
on the possibility and impossibility of BvM-theorems for undersmoothing
priors---his negative results will be relevant below. His positive
findings are, however, strongly tied to the Gaussian conjugate
situation, do not address efficiency questions, and do not give rise to
posteriors with the above mentioned ``plug-in property.'' For the related
question of obtaining semiparametric BvM-results, general sufficient
conditions are given in Castillo \cite{C12} and Bickel and Kleijn
\cite{BK12}, as well as in
Rivoirard and Rousseau \cite{RR12} for linear functionals of
probability density functions. A number of BvM-type results have been
obtained for the fixed finite-dimensional posterior with dimension
increasing to infinity: Ghosal \cite{G99} and Bontemps \cite{B11}
consider regression with a finite number of regressors, Ghosal \cite
{G00} and Clarke and Ghosal \cite{CG10} consider exponential families,
and the case of discrete probability distributions is treated in
Boucheron and Gassiat \cite{BG09}.

This article is organised as follows: in the next two subsections we
define a general notion of the nonparametric BvM phenomenon. In
Section \ref{confc} we demonstrate that when this phenomenon holds,
posterior-based inference is valid from a frequentist point of view in
a variety of concrete examples from nonparametric statistics. In
Section \ref{mainbvm} we prove that for a large class of natural
priors on $L^2$, the BvM phenomenon indeed occurs.

\subsection{The weak nonparametric Bernstein--von Mises phenomenon}

We consider a fixed design Gaussian regression model with known
variance, but work with its equivalent white noise formulation to
streamline the mathematical development. Let $L^2:=L^2([0,1])$ be the
space of square integrable functions on $[0,1]$. For $f \in L^2$, $dW$
standard white noise, consider observing
%
\begin{equation}
\label{wni} dX^{(n)}(t) = f(t)\,dt + \frac{1}{\sqrt{n}}\,dW(t),\qquad t\in[0,1].
\end{equation}
Except in conjugate situations the proof of a Bernstein--von Mises-type
result rests typically on the fact that efficient estimation at the
rate $1/\sqrt n$ is possible. In the nonparametric situation this rules
out $L^p$-type loss functions, but leads one to consider weaker $\ell
^\infty(\mathcal H)$-type norms discussed in (\ref{pipu}). For the
particular choice of $\mathcal H_s$ equal to an order-$s$ Sobolev-ball,
we can understand this better by using simple but useful Hilbert space
duality arguments in the nested scale of Sobolev spaces $\{H^r_2\}_{r
\in\mathbb R}$ on $[0,1]$: we define these in precise detail below,
but note for the moment that $H^r_2 \subseteq H^t_2, r \ge
t, H^0=L^2$, so to weaken the norm beyond $L^2$ means that we should
decrease $r$ to be negative. For $s>0$ the space $H^{-s}_2$ can be
realised in an isometric way as a closed subspace of $l^\infty
(\mathcal H_s)$, explaining heuristically the connection to the
discussion surrounding (\ref{pipu}) above. The space should be large
enough so that the Gaussian experiment in (\ref{wni}) can be realised
as a tight random element in $H^{-s}_2$. The critical value for this to
be the case is $s=1/2$, and we define in (\ref{spaces}) below a (in a
certain sense ``maximal'') Sobolev space $H$ with norm $\|\cdot\|_H$ in
which the random trajectory $dX^{(n)}$ defines a tight Gaussian Borel
random variable $\mathbb{X}^{(n)}$ with mean $f$ and covariance
$n^{-1}I$. That
is, if we denote by $\dob$ the centered Gaussian Borel random variable
on $H$ with covariance $I$, then (\ref{wni}) can be written as
%
\begin{equation}
\label{shift} \mathbb{X}^{(n)}= f + \frac{1}{\sqrt{n}}\dob,
\end{equation}
a natural Gaussian shift experiment in the Hilbert space $H$. One can
show moreover that $\mathbb{X}^{(n)}$ is an efficient estimator for
$f$ for the
loss function of $H$.

Any (Borel or cylindrical) probability measure on $L^2$ gives rise to a
tight probability measure on $H$ simply by the continuous
(Hilbert--Schmidt) injection $L^2 \subset H$. Let thus $\Pi$ be a
prior on $L^2$, and let
\[
\Pi_n=\Pi\bigl(\cdot|X^{(n)}\bigr)=\Pi\bigl(\cdot|
\mathbb{X}^{(n)}\bigr)
\]
be the posterior distribution on $H$ given the observed trajectory from
(\ref{wni}), or equivalently, from (\ref{shift}). On $H$ and for $z
\in H$, define the transformation
\[
\tau_z\dvtx  f \mapsto\sqrt n (f-z).
\]
Let $\Pi_n \circ\tau_{\mathbb{X}^{(n)}}^{-1}$ be the image of the
posterior law
under $\tau_{\mathbb{X}^{(n)}}$. The shape of $\Pi_n \circ\tau
_{\mathbb{X}^{(n)}}^{-1}$
reveals how the posterior concentrates on $1/\sqrt n$-$H$-neighborhoods
of the efficient estimator $\mathbb{X}^{(n)}$. To compare probability
distributions on $H$ we may use any metric for weak convergence of
probability measures, and we choose the bounded Lipschitz metric here
for convenience (it is defined in Section \ref{weak}). Let $\mathcal
N$ be the standard Gaussian probability measure on $H$ with mean zero
and covariance $I$ constructed in Section \ref{sob} below. It should
be distinguished from the standard Gaussian law $N(0,I)$ on $\mathbb
R^k, k \in\mathbb N$.

\begin{definition}\label{bvmp}
Consider data generated from equation (\ref{wni}) under a fixed
function $f_0$, and denote by $P_{f_0}^n$ the distribution of $\mathbb
{X}^{(n)}$.
Let $\beta$ be the bounded Lipschitz metric for weak convergence of
probability measures on $H$. We say that a prior $\Pi$ satisfies the
weak Bernstein--von Mises phenomenon in $H$ if, as $n \to\infty$,
\[
\beta\bigl(\Pi_n \circ\tau_{\mathbb{X}^{(n)}}^{-1},\mathcal N
\bigr) \to^{P_{f_0}^n}0.
\]
\end{definition}

We note that the fact that the result is phrased in a way in which
$\mathcal N$ is independent of $n$ is important since $\beta$ does not
induce a uniformity structure for the topology of weak convergence; see
the remark on page 413 in \cite{D02}.

Thus when the weak Bernstein--von Mises phenomenon holds, the posterior
necessarily has the approximate shape of an infinite-dimensional
Gaussian distribution. Moreover, we require this Gaussian distribution
to equal $\cN$, the canonical choice in view of efficiency
considerations. The covariance of $\cN$ is the Cram\'er--Rao bound for
estimating $f$ in the Gaussian shift experiment (\ref{shift}) in
$H$-loss, and we shall see how this carries over to sufficiently
regular real-valued functionals $\Psi(f)$; see Section \ref
{sec-funct} below.

One may ask by analogy to the finite-dimensional situation whether a
\textit{strong} Bernstein--von Mises phenomenon, where $\beta$ is
replaced by the total variation norm, can be proved. It follows from
Theorem 2 in \cite{L11} that already in the Gaussian conjugate
situation, such a result is impossible unless one restricts to very
specific priors (which in particular do not possess the plug-in
property that will be needed in the key applications below).

Now with weak instead of total variation convergence, we cannot infer
that $\Pi_n \circ\tau_{\mathbb{X}^{(n)}}^{-1}$ and $\mathcal N$ are
approximately\vspace*{1pt} the same for every Borel set in $H$, but only for sets
$B$ that are continuity sets for the probability measure $\mathcal N$.
For statistical applications of the Bernstein--von Mises phenomenon,
one typically needs some uniformity in $B$, and this is where total
variation results would be particularly useful. Weak convergence in $H$
implies that $\Pi_n \circ\tau_{\mathbb{X}^{(n)}}^{-1}$ is close to
$\mathcal N$
uniformly in certain classes of subsets of $H$ whose boundaries are
sufficiently regular relative to the measure~$\mathcal N$ (see
Section \ref{weak}), and we show below how this allows for enough
uniformity to deal with a variety of concrete nonparametric statistical
problems.

The Bernstein--von Mises phenomenon in Definition \ref{bvmp} will
often be complemented by convergence of moments, that is, convergence
of the Bochner integrals (e.g., page 100 in \cite{AG80})
$
\int_{H} f\,d\Pi_n \circ\tau_{\mathbb{X}^{(n)}}^{-1}(f)
\to^{P_{f_0}^n} \int_{H} f\,d\mathcal N(f) =0$
as $n \to\infty$ in~$H$. This implies that the posterior mean $\bar
f_n$ of $\Pi_n$ satisfies
%
\begin{equation}
\label{asylin} \bigl\|\bar f_n - \mathbb{X}^{(n)}
\bigr\|_{H} = o_P\bigl(n^{-1/2}\bigr),
\end{equation}
so in semiparametric terminology the posterior mean is asymptotically
linear in $H$ with respect to 
$\mathbb{X}^{(n)}$; in particular, $\bar f_n$ is an efficient
estimator for $f$.

\subsection{Sobolev spaces and white noise} \label{sob}

Denote by $\langle f, g \rangle= \int_0^1 f(x) \overline{g(x)}\,dx$
the standard inner product on $L^2$. General order Sobolev spaces will
be defined via orthonormal bases of $L^2$ that satisfy the following
weak regularity condition. While notationally it reflects a
wavelet-type basis $\{\psi_{lk}\dvtx  l \ge J_0-1, 0 \le k \le2^l-1\}$ of
CDV-type \cite{CDV93} (with notational convention $\psi_{(J_0-1)k} =
\phi_k$ for the scaling function), it also includes the trigonometric
basis $\psi_{lk}(x)\equiv e_l(x) = e^{2\pi i l x}$ and bases of
standard Karhunen--Lo\`eve expansions.

\begin{definition}\label{onb}
Let $S \in\mathbb N$. By an $S$-regular basis $\{\psi_{lk}\dvtx  l \in
\mathcal L, k \in\mathcal Z_l\}$ of $L^2$ with index sets $\mathcal L
\subset\mathbb Z, \mathcal Z_l \subset\mathbb Z$ and characteristic
sequence $a_l$ we shall mean any of the following:

(a) $\psi_{lk}\equiv e_{l}$ is $S$-times differentiable with
all derivatives in $L^2$, $|\mathcal Z_l|=1$, $a_l = \max(2, |l|)$,
and $\{e_l\dvtx  l \in\mathcal L\}$ forms an orthonormal basis of $L^2$.

(b) $\psi_{lk}$ is $S$-times differentiable with all
derivatives in $L^2$, $\mathcal L \subset\mathbb N$, $a_l=|\mathcal
Z_l|=2^l$, and $\{\psi_{lk}\dvtx  l \in\mathcal L, k \in\mathcal Z_l\}$
forms an orthonormal basis of $L^2$.
\end{definition}
Define for $0 \le s < S$ the standard Sobolev spaces as
\[
H^s_2:= \biggl\{f \in L^2\bigl([0,1]
\bigr)\dvtx  \|f\|^2_{s,2}:= \sum_{l \in
\mathcal L}
a_l^{2s} \sum_{k \in\mathcal Z_l} \bigl|\langle
\psi_{lk}, f \rangle\bigr|^2 < \infty\biggr\},
\]
which for the usual wavelet or trigonometric bases are in fact spaces
independent of the basis. For $\mathcal Z_l' \subset\mathcal Z_l,
\mathcal L' \subset\mathcal L$ finite we can form linear subspaces
\[
V \equiv V_{\mathcal L', \mathcal Z_l'}=\operatorname{span}\bigl\{\psi_{lk}\dvtx  l \in
\mathcal
L', k \in\mathcal Z_l'\bigr\}
\]
of $H^s_2 \subset L^2$, and we denote the $L^2$-projection of $f \in
L^2$ onto $V$ by $\pi_{V}(f)$.

For $s>0$ we define the dual space
\[
H^{-s}_2\bigl([0,1]\bigr):= \bigl(H^s_2[0,1]
\bigr)^*.
\]
Using standard duality arguments (as in Proposition 9.16 in \cite
{F99}) one shows the following: $H^{-s}_2$ consists precisely of those
linear forms $L$ acting on $H^s_2$ for which the $\|L\|_{-s,2}$-norms
[defined as above also for negative $s$, with $\langle\psi_{lk}, L
\rangle$ replaced by $L(\psi_{lk})$, noting $\psi_{lk} \in H^s_2$]
are finite. In fact the so-defined norm $\|\cdot\|_{-s,2}$ is
equivalent to the standard operator norm on $(H^s_2[0,1])^*$. Moreover
every $f \in L^2$ gives rise to a continuous linear form on $H^s_2
\subset L^2$ by using the $\langle\cdot, \cdot\rangle$ duality, so
we can view $L^2$ as a subspace of $H^{-s}_2$. By reflexivity of
$H^s_2$ one concludes $(H^{-s}_2([0,1]))^* = H^s_2([0,1])$ up to
isomorphism, that is, any linear continuous map $K\dvtx H^{-s}_2 \to\mathbb
R$ is of the form $K\dvtx  L \mapsto L(g)$ for some $g \in H^s_2$, and if
$L$ itself is a functional coming from integrating against an
$L^2$-function $f_L$, then $L(g)=\langle g, f_L \rangle$.

To obtain sharp results we also need ``logarithmic'' Sobolev spaces
\[
H^{s,\delta}_2 \equiv\biggl\{f\dvtx  \|f\|^2_{s,2,\delta}:= \sum_{l
\in\mathcal L} \frac{a^{2s}_l} {(\log a_l)^{2\delta}} \sum
_{k \in
\mathcal Z_l} \bigl|\langle\psi_{lk}, f \rangle\bigr|^2
< \infty\biggr\},\qquad \delta\ge0, s \in\mathbb R,
\]
which are Hilbert spaces satisfying the compact imbeddings $H^{r}_2
\subset H^{r,\delta}_2 \subset H^{s}_2$ for any real valued $s<r$.

For any $f \in H^s_2 \subseteq L^2$ ($s \ge0$) and $dW$ standard white
noise, we have a random linear application
%
\begin{equation}
\label{dual} \dob\dvtx  f \mapsto\int_0^1 f(t)
\,dW(t) \sim N\bigl(0, \|f\|_2^2\bigr).
\end{equation}
For any $\delta>1/2$, the $\|\dob\|_{-1/2,2,\delta}$-norm converges
almost surely since, by Fubini's theorem, for $g_{lk}$ independent
$N(0,1)$ variables,
\[
E \|\dob\|^2_{-1/2,2,\delta} = \sum_{l \in\mathcal L}
a_l^{-1} (\log a_l)^{-2\delta} \sum
_{k \in\mathcal Z_l} E{g_{lk}^2} < \infty,
\]
so $\dob\in H^{-1/2, \delta}_2$ almost surely, measurable for the
cylindrical $\sigma$-algebra, and by separability of $H^{-1/2,\delta
}_2$ also for the Borel $\sigma$-algebra (page 374 in \cite{B98}). By
Ulam's theorem (Theorem 7.1.4 in \cite{D02}), $\dob$ is thus tight in
$H^{-1/2,\delta}_2$. One can show that the spaces
%
\begin{equation}
\label{spaces}\qquad H\equiv H(\delta)\equiv H^{-1/2,\delta}_2,\qquad \|\cdot
\|_H \equiv\| \cdot\|_{H(\delta)} \equiv\|\cdot
\|_{-1/2,2,\delta},\qquad \delta>1/2,
\end{equation}
are minimal in the considered scale of spaces on which this happens:
decreasing $\delta$ below $1/2$ would lead to a space in which $\dob$
is not tight.

The Gaussian variable $\dob$ has mean zero and covariance $I$ diagonal
for the $L^2$-inner product, that is,
$E\dob(g)\dob(h) = \langle g, h \rangle$, for all $g,h \in L^2$. We
call the law $\mathcal N$ of $\dob$ a standard, or canonical, Gaussian
probability measure on the Hilbert space $H$ (note that it is the
isonormal Gaussian measure for the inner product of $L^2$ but \textit
{not} for the one of $H$). In the same way the random trajectory
$dX^{(n)}$ from (\ref{wni}) defines a tight Gaussian Borel random
variable $\mathbb{X}^{(n)}$ on $H$ with mean $f$ and covariance
$n^{-1}I$, thus
rigorously justifying (\ref{shift}).

We finally define H\"{o}lder-type spaces of smooth functions: for
$S>s>0$ and $\psi_{lk}$ a $S$-regular wavelet basis from Definition
\ref{onb}(b), we set
%
\begin{equation}
\label{hold} C^s \equiv\Bigl\{f \in C\bigl([0,1]\bigr)\dvtx  \|f
\|_{s, \infty}:= \sup_{l \in
\mathcal L, k \in\mathcal Z_l} 2^{l(s+1/2)} \bigl|\langle
\psi_{lk}, f \rangle\bigr| < \infty\Bigr\}.
\end{equation}

\section{Confidence sets for nonparametric Bayes procedures} \label{confc}

\subsection{Weighted $L^2$-credible ellipsoids}

Throughout this section $H$ stands for the space $H(\delta)$ from
(\ref{spaces}) for some arbitrary choice of $\delta>1/2$.
Denote by $B(g,r)= \{f \in H\dvtx  \|f-g\|_H \le r\}$ the norm ball in $H$
of radius $r$ centered at $g$. In terms of an orthonormal basis $\{\psi
_{lk}\}$ of $L^2$ from Definition \ref{onb} this corresponds to
$L^2$-ellipsoids
\[
\biggl\{\{c_{lk}\}\dvtx  \sum_{l,k}
a_l^{-1} \bigl(\log(a_l)\bigr)^{-2\delta}
\bigl|c_{lk}-\langle g, \psi_{lk} \rangle\bigr|^2 \le
r^2 \biggr\},
\]
where coefficients in the tail are downweighted by $a_l^{-1} (\log
(a_l))^{-2\delta}$. A frequentist goodness of fit test of a null
hypothesis $H_0\dvtx  f=f_0$ could, for instance, be based on the test
statistic $\|f_0-\mathbb{X}^{(n)}\|_H$, resembling in nature a Cram\'er--von
Mises-type procedure that has power against arbitrary fixed
alternatives $f \in L^2$.

A Bayesian approach consists in using the quantiles of the posterior
directly. Given $\alpha>0$ one solves for $R_n\equiv R ( \mathbb
{X}^{(n)}, \alpha
)$ such that
%
\begin{equation}
\label{qut} \Pi\bigl(f\dvtx  \| f - T_n \|_{H} \le
R_n/\sqrt n\given \mathbb{X}^{(n)}\bigr) = 1-\alpha,
\end{equation}
where $T_n=\mathbb{X}^{(n)}$ or, when the posterior mean $\bar f_n$ exists,
possibly $T_n=\bar f_n$. A~$\|\cdot\|_{H}$-ball centred at $T_n$ of
radius $R_n$ constitutes a level $(1-\alpha)$-credible set for the
posterior distribution. The weak Bernstein--von Mises phenomenon in $H$
implies that this credible ball asymptotically coincides with the exact
$(1-\al)$-confidence set built using the efficient estimator $\mathbb
{X}^{(n)}$
for $f$.

\begin{theorem} \label{sct}
Suppose the weak Bernstein--von Mises phenomenon in the sense of
Definition \ref{bvmp} holds. Given $0<\alpha<1$ consider the credible set
%
\begin{equation}
\label{cred} C_n= \bigl\{f\dvtx  \bigl\|f-\mathbb{X}^{(n)}
\bigr\|_{H} \le R_n/\sqrt n \bigr\},
\end{equation}
where $R_n\equiv R ( \mathbb{X}^{(n)}, \alpha)$ is such that $\Pi
(C_n| \mathbb{X}^{(n)}) =
1-\alpha$. Then, as $n \to\infty$,
\[
P^n_{f_0} (f_0 \in C_n) \to1-
\alpha\quad\mbox{and}\quad R_n=O_P(1).
\]
If in addition $\|\bar f_n - \mathbb{X}^{(n)}\|_{H} = o_P(n^{-1/2})$,
then the
same is true if in the definition of $C_n$ the posterior mean $\bar
f_n$ replaces $\mathbb{X}^{(n)}$.
\end{theorem}

When available, using further prior knowledge in the construction of
the credible set may lead to favourable frequentist properties, such as
optimal performance in stronger loss functions.

To see this, consider first the specific but instructive case of a
uniform wavelet prior $\Pi$ on $L^2$ arising from the law of the
random wavelet series
\[
U_{\gamma,M} = \sum_{l = J_0-1}^\infty\sum
_{k =0}^{2^l-1} 2^{-l(\gamma+1/2)}u_{lk}
\psi_{lk}(\cdot),\qquad \gamma>0,
\]
where the $u_{lk}$ are i.i.d. uniform on $[-M,M]$ for some $M>0$, with
$S$-regular CDV-wavelets $\psi_{lk}$, $S>\max(\gamma,1/2), J_0 \in
\mathbb N$. Such priors model functions that lie in a fixed H\"older
ball of $\|\cdot\|_{\gamma, \infty}$-radius $M$, with posteriors
$\Pi(\cdot|\mathbb{X}^{(n)})$ contracting about $f_0$ at the
$L^2$-minimax rate
within logarithmic factors if $\|f_0\|_{\gamma, \infty} \le M$; see
\cite{GN11} and also Corollary \ref{postrat} below.

In this situation it is natural to intersect the credible set $C_n$
with the H\"olderian support of the prior (or posterior),
%
\begin{equation}
\label{credunif} C'_n = \bigl\{f\dvtx  \|f\|_{\gamma, \infty} \le
M, \|f-\bar f_n\|_{H} \le R_n/\sqrt n \bigr\},
\end{equation}
where $R_n$ is as in (\ref{qut}) with $T_n= \bar f_n$. Note that the
posterior mean also satisfies $\|\bar f_n\|_{\gamma, \infty} \le M$,
so that $C_n'$ is a random subset of a H\"older ball that has
credibility $\Pi(C'_n|\mathbb{X}^{(n)})=1-\alpha$. Theorem \ref
{sct} implies
the following result.

\begin{corollary} \label{unifset}
Consider observations generated from (\ref{shift}) under a fixed
function $f_0 \in C^\gamma$ with $\|f_0\|_{\gamma, \infty}<M$. Let
$\Pi$ be the law of $U_{\gamma, M}$, let $\Pi(\cdot|\mathbb
{X}^{(n)})$ be the
posterior distribution given $\mathbb{X}^{(n)}$ and let $C'_n$ be as
in (\ref
{credunif}). Then
\[
P^n_{f_0} \bigl(f_0 \in C'_n
\bigr) \to1-\alpha
\]
as $n \to\infty$ and the $L^2$-diameter $|C'_n|_2$ of $C'_n$
satisfies, for some $\kappa>0$,
\[
\bigl|C'_n\bigr|_2 = O_P
\bigl(n^{-\gamma/(2\gamma+1)} (\log n) ^\kappa\bigr).
\]
\end{corollary}

We consider next the situation of a general series prior $\Pi$
modelling $\gamma$-regular functions, including the important case of
Gaussian priors. Let
\[
G_{\gamma} = \sum_{l = J_0-1}^\infty\sum
_{k =0}^{2^l-1} 2^{-l(\gamma+1/2)}
g_{lk} \psi_{lk}(\cdot),\qquad \gamma>0,
\]
where $g_{lk}$ are i.i.d. random variables that possess a bounded
positive density $\varphi$ such that $\operatorname{Var}(g_{lk})<\infty$, and with
$S$-regular CDV-wavelets $\psi_{lk}$, $S>\max(\gamma,1/2)$. Denote
by $\Pi_n = \Pi(\cdot|\mathbb{X}^{(n)})$ the posterior distribution from
observing $\mathbb{X}^{(n)}\sim P_{f_0}^n$. The idea behind the
construction of
$C_n'$ can be adapted to this more general situation by taking for $M_n
\to\infty,  M_n=O(\log n)$,
\[
\tilde C_n' = \bigl\{f\dvtx  \|f\|_{\gamma, 2} \le
M_n, \|f-\bar f_n\|_{H} \le R_n/
\sqrt n \bigr\}.
\]
This parallels the frequentist practice of ``undersmoothing,'' taking
into account the fact that we usually do not know a bound for $\|f\|
_{\gamma, 2}$ in the construction of confidence sets. While this can
be shown to work (as in the proof of Corollary \ref{unifset}, assuming
$f_0 \in C^\gamma\cap H^\gamma_2$), we wish to avoid such ad hoc
methods here and prefer to explicitly use posterior information on the
size of $\|f\|_{\gamma, 2,1}$: fix $\delta>0$ arbitrarily, and set
%
\begin{equation}
\label{credgauss} C_n''= \bigl\{ f\dvtx  \|f
\|_{\gamma,2,1} \le M_n+4\delta, \|\bar{f}_n - f\|
_H \le R_n / \sqrt{n} \bigr\},
\end{equation}
where $R_n$ is as in (\ref{qut}) with $T_n=\bar f_n$, and where $M_n$
is defined as follows: for any $n$ and $\delta_n=(\log{n})^{-1/4}$,
%
\begin{equation}
\label{defqn2} M_n = \inf\bigl\{ M>0\dvtx   \Pi_n\bigl(f\dvtx  \bigl|
\|f\|_{\gamma,2,1} - M \bigr|\le\delta\bigr) \ge1-\delta_n \bigr\}
\end{equation}
with the convention that $M_n=\infty$ if the set over
which one takes the infimum in (\ref{defqn2}) is empty.


\begin{corollary} \label{gaussset}
Consider observations generated from equation (\ref{shift}) under a
fixed function
$f_0 \in C^\gamma$. Let $\Pi$ be the law of $G_{\gamma}$, let $\Pi
_n=\Pi(\cdot|\mathbb{X}^{(n)})$ be the posterior distribution given
$\mathbb{X}^{(n)}$, and
let $C_n''$ be as in (\ref{credgauss}). Then
\[
P^n_{f_0} \bigl(f_0 \in C''_n
\bigr) \to1-\alpha,\qquad \Pi_n\bigl(C''_n
\bigr) = 1-\alpha+ o_P(1)
\]
as $n \to\infty$, and the $L^2$-diameter $|C''_n|_2$ of $C''_n$
satisfies, for some $\kappa>0$,
\[
\bigl|C''_n\bigr|_2 = O_P
\bigl(n^{-\gamma/(2\gamma+1)} (\log n) ^\kappa\bigr).
\]
Additionally, both $M_n$ and $R_n$ occuring in (\ref{credgauss}) are
bounded in probability.
\end{corollary}

These credible sets can be compared to those in \cite{KVV11}
constructed in the Gaussian conjugate situation [i.e., for $g_{lk}$
i.i.d. $N(0,1)$]. Both constructions give rise to credible sets that
have frequentist minimax optimal diameter (within log-factors). In
contrast to $C_n''$, however, the credible sets in \cite{KVV11} are
conservative in the sense that their asymptotic frequentist coverage
probability may exceed the desired level $1-\alpha$.

The purpose of $M_n$ in (\ref{defqn2}) is to provide a bound on the
unknown $\|f\|_{\ga,2,1}$ using the posterior distribution, similar in
spirit\vspace*{1pt} to a posterior quantile. Using (\ref{defqn2}) and Theorem \ref
{thm-l2} below (with $\sigma_l= 2^{-l(\gamma+1/2)}$), one shows that
$\Pi_n(f\dvtx \|f-f_0\|_{\gamma,2,1} > (\log n)^{-1/4}) = o_P(1)$
and then also that
%
\begin{equation}
\label{normest} \|f_0\|_{\ga,2,1} -2\delta+
o_P(1) \le M_n \le\|f_0\|_{\ga,2,1} +
2\delta+o_P(1).
\end{equation}
It is also possible to take $\delta=\delta_n$ in (\ref{defqn2}).
Corollary \ref{gaussset} then still holds, and $\delta$ is replaced
by $\delta_n$ in the previous display, in which case $M_n$ is a
consistent estimator of $\|f\|_{\ga,2,1}$.

\subsection{Credible bands for self-convolutions}

We proceed with a semiparametric example: suppose we are interested in
estimating the function
\[
f \ast f = \int_0^1 f(\cdot-t)f(t)\,dt,
\]
where addition is $mod$-1 (so the convolution of $f$ with itself on the
unit circle). The related problem in density estimation was studied in
the papers \cite{F94,SW04,GM07,N07,N09}, where it is shown that $f
\ast f$ can be estimated at the $1/\sqrt n$-rate even when this is
impossible for $f$. See particularly \cite{F94} for applications.
Assume $f$ is one-periodic and contained in $H^{s}_2$ for some $s>1/2$,
and that the posterior is supported in $L^2([0,1))\equiv
L^2_{per}([0,1))$ which, in this subsection, denotes the subspace of
$L^2$ consisting of one-periodic functions. We will assume that the
basis used to define $H^s_2$ is such that $(\sum_{m}|\hat f(m)|^2
(1+|m|)^{2s})^{1/2}$ is an equivalent norm on $H^s_2$ (which is the
case for CDV- or periodised wavelets and trigonometric bases of $L^2$).

By standard properties of convolutions $\kappa\dvtx  f \mapsto f \ast f$
maps $L^2([0,1))$ into $C([0,1))$, the space of bounded continuous
periodic functions\vspace*{1pt} on $[0,1)$ equipped with the uniform norm $\|\cdot\|
_\infty$. If $\Pi_n = \Pi(\cdot|\mathbb{X}^{(n)})$ with posterior
mean $\bar
f_n \in L^2([0,1))$, we can construct a confidence band for $f \ast f$ by
solving for $R_n$ such that
%
\begin{equation}
\label{conv} \Pi_n \circ\kappa^{-1}\bigl(g\dvtx  \|g-\bar
f_n \ast\bar f_n\|_\infty\le R_n/
\sqrt n\bigr) =1-\alpha
\end{equation}
with resulting credible band
%
\begin{equation}
\label{cconv} C_n = \bigl\{g\dvtx  \|g-\bar f_n \ast\bar
f_n \|_\infty\le R_n/\sqrt n\bigr\}.
\end{equation}

\begin{theorem} \label{convt} Suppose the weak Bernstein--von Mises
phenomenon in the sense of Definition \ref{bvmp} holds, and that $f_0
\in H^s_2$ for some $s>1/2$. Assume
\[
\bigl\|\bar f_n - \mathbb{X}^{(n)}\bigr\|_{H} =
o_P\bigl(n^{-1/2}\bigr),
\]
and that for some sequence $r_n=o(n^{-1/2})$,
\[
\|\bar f_n - f_0\|_2^2 =
O_P(r_n),\qquad \Pi_n\bigl(f\dvtx \|f-f_0
\|_2^2 > r_n\bigr) =o_P(1).
\]
Let $C_n$ be the credible band from (\ref{cconv}) with $R_n$ as in
(\ref{conv}). Then, as $n \to\infty$,
\[
P^n_{f_0} (f_0 \ast f_0 \in
C_n) \to1-\alpha\quad\mbox{and}\quad R_n=O_P(1).
\]
\end{theorem}

If $f_0 \in C^s \cap H^s_2$ for some $s>1/2$, the priors from Condition
\ref{tolk} below with $\sigma_l$ and $\gamma=s$ chosen as in Remark
\ref{rates} are admissible in Theorem \ref{convt} with $r_n =
n^{-2s/(2s+1)}$; cf. Corollaries \ref{postrat}, \ref{barrat} below
and Section \ref{pipp}.


\subsection{Credible sets for functionals} \label{sec-funct}

\subsubsection{Linear functionals}

Let $L$ be any linear form on $L^2$ given by
\[
L(f)=\langle f, g_L\rangle= \int_0^1
f(t) g_L(t)\,dt,\qquad f \in L^2,
\]
where $g_L \in H^s_2, s>1/2$, and $g_L\neq0$. If $\Pi_n = \Pi(\cdot
|\mathbb{X}^{(n)})$ is the posterior, one may construct credible sets
for $L(f_0)$
based on the induced law $\Pi_n^L=\Pi_n \circ L^{-1}$ in several
ways: for example, one solves for $R_n=R(\mathbb{X}^{(n)}, L, \alpha
)$ in
%
\begin{equation}
\label{qtl} \Pi_n^L\bigl(z\dvtx  \bigl|z-L\bigl(
\mathbb{X}^{(n)}\bigr)\bigr| \le R_n/\sqrt n\bigr) = 1-\alpha,
\end{equation}
which gives rise to the credible set
%
\begin{equation}
\label{cfl} C_n = \bigl\{z\dvtx  \bigl|z-L\bigl(\mathbb{X}^{(n)}
\bigr)\bigr| \le R_n/\sqrt n \bigr\}
\end{equation}
for $L(f)$. An alternative way to build the credible set is discussed
below in a more general setting.
%
\begin{theorem} \label{linbvm}
Suppose the weak Bernstein--von Mises phenomenon in the sense of
Definition \ref{bvmp} holds. Let $L=\langle\cdot, g_L\rangle$ be a
linear functional on $L^2$ where $0 \neq g_L \in H^s_2, s>1/2$. Let
$\beta_{\mathbb{R}}$ be the bounded-Lipschitz metric for weak
convergence on $\mathbb R$, and define $\te_t\dvtx  x \mapsto\sqrt n
(x-t)$ for $t,x \in\mathbb R$. Then
\[
\be_{\mathbb{R}}\bigl( \Pi_n^L \circ
\te_{L(\mathbb{X}^{(n)})}^{-1}, N\bigl(0, \|g_L\|
_2^2\bigr) \bigr) \to^{P_{f_0}^n} 0.
\]
Moreover let $C_n$ be the credible region from (\ref{cfl}) with $R_n$
chosen as in (\ref{qtl}). Then
\[
P^n_{f_0} \bigl(L(f_0) \in C_n
\bigr) \to1-\alpha\quad\mbox{and}\quad R_n=O_P(1)
\]
as $n \to\infty$. If, in addition, (\ref{asylin}) holds, then the
same result holds true if $C_n$ is centered at $L(\bar f_n)$ where
$\bar f_n$ is the posterior mean of $\Pi(\cdot|\mathbb{X}^{(n)})$.
\end{theorem}

The induced posterior $\Pi_n\circ L^{-1}$ has the approximate shape of
a normal distribution centered at the efficient estimator $L(\mathbb
{X}^{(n)})$ of
$L(f)$ with variance $\|g_L\|_2^2/n$. This implies in particular that
the width of the credible set $C_n$ is asymptotically efficient from
the semiparametric perspective; in fact $\|g_L\|_2^2$ is the
semiparametric Cram\'er--Rao bound for estimating $L(f)$ from
observations in the Gaussian white noise model (when maintaining
standard nonparametric models for $f$).

The fact that any integral functional $\int f(t) g_L(t)\,dt, g_L \in
H^s_2, s>1/2$, is covered gives rise to a rich class of examples, such
as the moment functionals $\int t^\alpha f(t)\,dt, \alpha\in\mathbb
N$. The restriction to $s>1/2$ is intrinsic to our methods and cannot
be relaxed.

\subsubsection{Smooth nonlinear functionals}

We next consider statistical inference for nonlinear functionals of
$f_0$ that satisfy a good quadratic approximation in $L^2$ at $f_0$,
and more precisely, we assume that $\Psi\dvtx  L^2 \to\mathbb R$ satisfies
%
\begin{equation}
\label{qfuz} \Psi(f_0+h)-\Psi(f_0) = D
\Psi_{f_0}[h] + O\bigl(\|h\|_2^2\bigr),
\end{equation}
uniformly in $h \in L^2$ and for some $D\Psi_{f_0}\dvtx  L^2 \to\mathbb R$
linear and continuous that has a (nonzero) $L^2$-Riesz representer
$\dot{\Psi}_{f_0} \in H^s_2$ for some $s>1/2$. This setting includes
several standard examples discussed in more detail at the end of this
section, but also the linear functionals discussed above.

Note that now $\Psi$ cannot necessarily be evaluated at $\mathbb
{X}^{(n)}$ [think
of $\Psi(f) = \|f\|_2^2$].
However, since the posterior is supported in $L^2$ with probability
one, the following Bayesian credible set can be constructed for $\Psi
(f)$: for $\Pi_n = \Pi(\cdot|\mathbb{X}^{(n)})$ the posterior
distribution, set
$\Pi_n^\Psi=\Pi_n\circ\Psi^{-1}$, and solve for the $\al/2$ and
$1-\al/2$ quantiles $\mu_n,\nu_n$ of $\Pi_n^\Psi$,
%
\begin{equation}
\label{functz} \Pi_n^{\Psi}\bigl((-\infty, \mu_n] \bigr) =
\Pi_n^{\Psi}\bigl((\nu_n,+\infty) \bigr) =
\frac{\alpha}{2}.
\end{equation}

\begin{theorem} \label{bvm-func}
Suppose the weak Bernstein--von Mises phenomenon in the sense of
Definition \ref{bvmp} holds. Consider a functional $\Psi$ satisfying
(\ref{qfuz}). Assume moreover either that $\Psi$ is linear, or that
for some sequence $r_n=o(n^{-1/2})$, 
%
\[
\Pi_n\bigl(f\dvtx \|f-f_0\|_2^2 >
r_n\bigr) =o_P(1).
\]
Let $\mu_n,\nu_n$ satisfy (\ref{functz}). Then as $n \to\infty$,
\[
P^n_{f_0} \bigl(\Psi(f_0) \in(\mu_n,
\nu_n]\bigr) \to1-\alpha.
\]
\end{theorem}

Similar to Theorem \ref{linbvm}, the shape of the induced posterior
$\Pi_n \circ\Psi^{-1}$ is approximately Gaussian, this time centered
at $\Psi(f_0)+ \langle\dot{\Psi}_{f_0}/\sqrt n, \dob\rangle$,
and with variance $\|\dot{\Psi}_{f_0}\|_2^2/n$. More precisely, for
$\be_\RR$ the bounded-Lipschitz metric for weak convergence,
\[
\be_\RR\bigl(\Pi_n^{\Psi}\circ
\te^{-1}_{ \Psi(f_0)+
\langle\dot{\Psi}_{f_0}, \dob\rangle/{\rn}}, N\bigl(0, \|\dot{\Psi}_{f_0}
\|_2^2\bigr) \bigr) \to^{P_{f_0}^n}0.
\]
In fact\vspace*{1pt} the proofs imply that the random quantile $\mu_n$ admits the
expansion, for $\Phi_*$ the distribution function of a $N(0,\|\dot
{\Psi}_{f_0}\|^2)$ variable,
\[
\mu_n = \Psi(f_0)+ \frac{1}{\rn} \langle\dot{
\Psi}_{f_0}, \dob\rangle+ \frac{\Phi_*^{-1}(\al/2)}{\rn} +
o_P(1/\rn)
\]
and $\nu_n$ likewise, with $\Phi_*^{-1}(1-\frac\al2)$ replacing
$\Phi_*^{-1}(\frac\al2)$. Again, $\|\dot{\Psi}_{f}\|_2^2$ is the
semiparametric efficiency bound for estimating $\Psi(f)$ in the
Gaussian white noise model, which shows that the asymptotic width of
the credible set $(\mu_n, \nu_n]$ for $\Psi(f)$ is optimal in the
semiparametric sense.

If $f_0 \in C^\gamma$ for some $\ga>1/2$, then the priors from
Condition \ref{tolk} below with $\sigma_l$ chosen as in Remark \ref
{rates} are admissible in the above theorem with $r_n = n^{-2\ga/(2\ga
+1)}$; cf. Corollaries \ref{postrat}, \ref{barrat} and Section \ref
{pipp} below.

Examples include the standard quadratic functionals such as $\Psi
(f)=\int f^2(t)\,dt$ or composite functionals of the form $\Psi(f)=\int
\phi(f(x),x)\,dx$. Some functionals may necessitate some straightforward
modifications of our proofs: for instance, $\|f\|_p^p$ requires
differentiation on $L^p$ instead of $L^2$, and for the entropy
functional $\int f(t) \log f(t) \,dt$ one assumes $f_0 \ge\zeta>0$ on
$[0,1]$ and differentiates $\Psi$ on $L^\infty$. In these situations,
to control remainder terms, one may use contraction results in $L^p, 2
< p \le\infty$, instead of $L^2$, such as the ones in \cite{GN11}.
Our assumption $\ga>1/2$ is stronger than the critical assumption $\ga
\ge1/4$ needed for $1/\sqrt n$-estimability of some of these
functionals \cite{L96}, a~phenomenon intrinsic to general plug-in procedures.

\section{Bernstein--von Mises theorems in white noise} \label{mainbvm}

We now develop general tools that allow us to prove that priors satisfy
the Bernstein--von Mises phenomenon in the sense of Definition \ref
{bvmp}, and show how they can be successfully applied to a wide variety
of natural classes of product priors.

For $f \in L^2$ consider again observing a random trajectory in the
white noise model (\ref{wni}) of law $P_f^n$, with corresponding
expectation operator denoted by $E_f^n$. Given an orthonormal basis
from Definition \ref{onb}, the white noise model is equivalent to
observing the action of $\mathbb{X}^{(n)}$ on the basis, that is,
\[
\mathbb{X}^{(n)}_{lk} = \tlk+ \frac{1}{\rn}
\veps_{lk},\qquad k \in\mathcal Z_l, l \in\mathcal L,
\]
where $\theta_{lk} = \langle f, \psi_{lk}\rangle$, $\veps_{lk} \sim
^{\mathrm{i.i.d.}} N(0,1)$. Let $\Pi$ be a prior Borel probability distribution
on $L^2$ which induces a prior, also denoted by $\Pi$, on infinite
sequences $\{\te_{lk}\} \in l^2$. Let $\Pi(\cdot|\mathbb{X}^{(n)})$
be the
posterior distribution, and let $\Pi(\tlk\given\mathbb{X}^{(n)})$
denote the
marginal posterior on the coordinate $\tlk$.

\subsection{\texorpdfstring{Contraction results in $H(\delta)$}{Contraction results in H(delta)}}

In this subsection we consider priors of the form $\Pi=\bigotimes_{lk}
\pi_{lk}$ defined on the coordinates of the orthonormal basis $\{\psi
_{lk}\}$, where $\pi_{lk}$ are probability distributions with Lebesgue
density $\vphi_{lk}$ on the real line. Further assume, for some fixed
density $\vphi$ on the real line,
\[
\vphi_{lk}(\cdot) = \frac{1}{\sil} \vphi\biggl(\frac{\cdot}{\sil
}
\biggr)\qquad\forall k \in\mathcal Z_l \mbox{, with }
\sigma_l>0, \sum_{l,k} \sigma_l^2 <\infty.
\]

\begin{condition} \label{tolk}
(P1) Suppose that for a finite constant $M>0$,
\[
\sup_{l\in\mathcal L,k\in\mathcal Z_l} \frac{|\tolk|}{\sil} \le M.
\]

(P2) Suppose that $\vphi$ is such that for some $\ta>M$ and
$0<c_\vphi\le C_\vphi<\infty$
\begin{eqnarray*}
&\displaystyle \vphi(x)\le C_\vphi\qquad\forall x\in\RR,\qquad \vphi(x)\ge c_\vphi
\qquad\forall x\in(-\ta,\ta),&\\
&\displaystyle \int_\mathbb R x^2
\varphi(x)\,dx < \infty.&
\end{eqnarray*}
\end{condition}

Some discussion of this condition is in order: we allow for a rich
variety of base priors $\vphi$, such as Gaussian, sub-Gaussian,
Laplace, most Student laws, or more generally any law with positive
continuous density and finite second moment, but also uniform priors
with large enough support. The full prior on $f$ considered here is
thus a sum of independent terms over the basis $\{\psi_{lk}\}$,
including many, especially non-Gaussian, processes.
For Gaussian processes Condition \ref{tolk} applies simply by
verifying that the $L^2$-basis provided by the Karhunen--Lo\`eve
expansion of the process satisfies the conditions of Definition \ref
{onb}. This includes in particular Brownian motion: the corresponding
$\vphi$ is then the standard Gaussian density, and $\sil=1/(\pi
(l+\frac{1}{2}))$ are the square-roots of the eigenvalues of the
covariance operator. Through condition (P1), this allows for
signals $f_0\equiv(\te_{0,lk})$ whose coefficients on the basis
decrease at least as fast as $1/l$. For primitives of Brownian motion
similar remarks apply, with stronger but natural decay restrictions on
$\langle f_0, \psi_{lk} \rangle$.

In principle, making the prior rougher allows for more signals through
condition (P1), but this may harm the performance of the
posterior in stronger loss functions than the one considered in the
next theorem. Its proof basically consists of showing that, under
$P_{f_0}^n$, the coordinate-wise marginal posterior distributions
contract about each ``true'' coordinate $\langle f_0, \psi_{lk} \rangle
$ at rate $1/\sqrt n$ with constants independent of $k,l$.

\begin{theorem} \label{thm-mins}
Consider data generated from equation (\ref{wni}) under a fixed
function $f_0 \in L^2$ with coefficients
$\theta_0=\{\tolk\}=\{\langle f_0, \psi_{lk} \rangle\}$. Then if
the product prior $\Pi$ and $f_0$ satisfy Condition \ref{tolk}, we
have for every $\delta>1/2$, as $n\to\infty$,
\[
E^{n}_{f_0} \int\|f-f_0\|_{H(\delta)}^2
\,d\Pi\bigl(f \given\mathbb{X}^{(n)}\bigr) = O \biggl(\frac{1}{n}
\biggr).
\]
\end{theorem}
\begin{pf}
We decompose the index set $\mathcal L$ into $\cJ_n:=\{l\in
\mathcal L, \rn\sigma_l \ge S_0 \}$ and its complement, where $S_0$
is a fixed positive constant. The quantity we wish to bound equals, by
definition of the $H$-norm and Fubini's theorem,
\[
\sum_{l,k} a_l^{-1} (\log
a_l)^{-2\delta} E_{f_0}^{n} \int(\tlk-
\tolk)^2 \,d\Pi\bigl(\tlk\given\mathbb{X}^{(n)}\bigr).
\]
Define further
$B_{lk}(\mathbb{X}^{(n)}):= \int(\tlk-\tolk)^2 \,d\Pi(\tlk\given
\mathbb{X}^{(n)})$ whose
$P_{f_0}^n$-expectation we now bound. We write $\ix=\mathbb{X}^{(n)}$ and
$E=E_{f_0}^{n}$ throughout the proof to ease notation.

Using the independence structure of the prior we have $\Pi(\tlk\given
\ix)=\break\pi_{lk}(\tlk\given\ix_{lk})$, and under $P_{f_0}^{n}$,
\begin{eqnarray*}
B_{lk}(\ix) &=& \frac{\int(\tlk- \tolk)^2
e^{-{n}(\tlk- \tolk)^2/{2} +\rn\elk(\tlk- \tolk)} \vphi
_{lk}(\tlk)\,d\tlk} {
\int e^{-{n}(\tlk- \tolk)^2/{2} +\rn\elk(\tlk- \tolk)}
\vphi_{lk}(\tlk)\,d\tlk}
\\
&=& \frac{1}{n} \frac{\int v^2
e^{-{v^2}/{2} +\elk v} \vphi(
({\tolk+n^{-1/2}v})/{\sil} )/({\rn\sil})\,dv} {
\int e^{-{v^2}/{2} +\elk v} \vphi
(({\tolk+n^{-1/2}v})/{\sil} )/({\rn\sil})\,dv} =:\frac{1}{n}\frac{N_{lk}}{D_{lk}}(
\elk).
\end{eqnarray*}
About the indices $l \in \mathcal J_n^c$:
Taking a smaller integrating set on the denominator makes the integral smaller
\[
D_{kl}(\veps_{lk}) \ge\int_{-\sqrt{n} \sil}^{\sqrt{n} \sil}
e^{-{v^2}/{2} +\elk
v}\frac{1}{\rn\sil} \vphi\biggl(\frac{\tolk+n^{-1/2}v}{\sil
} \biggr) \,dv.
\]
To simplify the notation we suppose that $\ta>M+1$. If this is not the
case, one multiplies the bounds of the integral in the last display by
a small enough constant.
The argument of the function $\vphi$ in the
previous display stays in $[-M+1,M+1]$ under (P1). Under
assumption (P2) this implies that the value of $\vphi$ in the
last expression is bounded from below by $c_\vphi$.
Next applying Jensen's inequality with the logarithm function
%
\begin{eqnarray*}
\log D_{kl}(\veps_{lk}) &\ge& \log(2c_\vphi)-
\int_{-\sqrt{n} \sil}^{\sqrt{n} \sil} \frac{v^2}{ 2}\,
\frac{dv}{2\sqrt{n}\sil} + \elk\int_{-\sqrt{n} \sil}^{\sqrt{n} \sil}
v\,\frac
{dv}{2\sqrt{n}\sil}
\\
&=& \log(2c_\vphi)- (\sqrt{n}\sil)^2/6.
\end{eqnarray*}
Thus, $D_{kl}(\veps_{lk}) \ge2c_\vphi e^{- (\sqrt{n}\sil)^2/6}$,
which is bounded away from zero for indices in~$\cJ_n^c$.
Now about the numerator, let us split the integral defining $N_{kl}$
into two parts
$\{v\dvtx  |v|\le\rn\sil\}$ and $\{v\dvtx  |v| > \rn\sil\}$. That is,
$N_{kl}(\veps_{lk})=(I)+(\mathit{II})$. Taking the expectation of the first
term and using Fubini's theorem,
\[
E (I) = \int_{-\sqrt{n} \sil}^{\sqrt{n} \sil} v^2
e^{-{v^2}/{2}} E\bigl[ e^{\elk v}\bigr] \frac{1}{\rn\sil}\vphi
\biggl(
\frac{\tolk+n^{-1/2}v}{\sil} \biggr)\,dv \le2 n\sil^2C_\vphi/3.
\]
The expectation of the second term is bounded by first applying
Fubini's theorem as before and then changing variables back,
\begin{eqnarray*}
E (\mathit{II}) &=& \int_{|v|> \sqrt{n} \sil} v^2 e^{-{v^2}/{2}} E
\bigl[ e^{\elk v}\bigr] \frac{1}{\rn\sil}\vphi\biggl(\frac{\tolk
+n^{-1/2}v}{\sil}
\biggr)\,dv
\\
&=& \int_{{\tolk}/{\sil}+1}^{\pli} \biggl(\rn\sil u - \rn\sil
\frac{\tolk}{\sil} \biggr)^2\vphi(u)\,du
\\
&&{} + \int_{-\infty}^{{\tolk}/{\sil}-1} \biggl(\rn\sil u - \rn\sil
\frac{\tolk}{\sil} \biggr)^2\vphi(u)\,du
\\
&\le& 2 n\sil^2 \biggl[ \frac{\tolk^2}{\sil^2} + \int
_{-\infty
}^{\pli} u^2\vphi(u)\,du \biggr].
\end{eqnarray*}
Thus, using (P1) again, $E (I)+E(\mathit{II})$ is bounded on ${\cJ^c_n}$
by a fixed constant times~$n\sigma_l^2$. In particular, there exists a
fixed constant independent of $n, k, l$ such that $E(nB_{lk}(X))$ is
bounded from above by a constant on ${\cJ^c_n}$.

Now about the indices in $\cJ_n$. For such $l,k$, using
(P1)--(P2), one can find $L_0>0$ depending only on $S_0, M, \ta$ such
that, for any $v$ in $(-L_0,L_0)$, $\vphi((\tolk+n^{-1/2}v)/\sil)
\ge
c_\vphi$. Thus the denominator $D_{lk}(\elk)$ can be bounded from
below by
\[
D_{lk}(\elk) \ge c_\vphi\int_{-L_0}^{L_0}
e^{-{v^2}/{2} +\elk v}\frac{1}{\rn\sil} \,dv.
\]
On the other hand, the numerator can be bounded above by
\[
N_{lk}(\elk) \le C_\vphi\int v^2
e^{-{v^2}/{2} +\elk v}\frac{1}{\rn\sil} \,dv.
\]
Putting these two bounds together leads to
\[
B_{lk}(\elk) \le\frac{1}{n} \frac{ C_\vphi}{ c_\vphi}
\frac{ \int v^2 e^{-{v^2}/{2} +\elk v}\,dv}{\int_{-L_0}^{L_0}
e^{-{v^2}/{2} +\elk v} \,dv}.
\]
The last quantity has a distribution independent of $l,k$. Let us thus
show that
\[
Q(L_0)= E \biggl[ \frac{ \int v^2 e^{-(v-\veps)^2/2}\,dv}{\int
_{-L_0}^{L_0}
e^{-(v-\veps)^2/2} \,dv} \biggr]
\]
is finite for every $L_0>0$, where $\veps\sim N(0,1)$.
In the numerator we substitute $u=v-\veps$. Using the inequality
$(u+\elk)^2\le2v^2+2\elk^2$, the second moment of a standard normal
variable appears, and this leads to the bound
\[
Q(L_0) \le C E \biggl[ \frac{1+\veps^2}{ \int_{-L_0}^{L_0}
e^{-(v-\veps)^2/2} \,dv} \biggr]
\]
for some finite constant $C>0$. Denote by $g$ the density of a standard
normal variable, by $\Phi$ its distribution function and $\bar{\Phi
}=1-\Phi$. It is enough to prove that the following quantity is finite:
\begin{eqnarray*}
q(L_0)&:=& \int_{-\infty}^{+\infty}
\frac{(1+u^2)g(u)}{\bfi
(u-L_0)-\bfi(u+L_0)}\,du \\
&=& 2 \int_{0}^{+\infty}
\frac{(1+u^2)g(u)}{\bfi(u-L_0)-\bfi
(u+L_0)}\,du,
\end{eqnarray*}
since the integrand is an even function. Using the standard inequalities
\[
\frac{1}{\sqrt{2\pi}} \frac{u^2}{1+u^2}\frac{1}{u} e^{-u^2/2} \le
\bfi(u) \le\frac{1}{\sqrt{2\pi}} \frac{1}{u} e^{-u^2/2},\qquad u\ge1,
\]
it follows that for any $\delta>0$, one can find $M_\delta>0$ such
that, for any $u\ge M_\delta$, 
\[
(1-\delta) \frac{1}{u} e^{-u^2/2} \le\sqrt{2\pi} \bfi(u) \le
\frac{1}{u} e^{-u^2/2},\qquad u\ge M_\delta.
\]
Set $A_\delta=2L_0 \vee M_\delta$. Then for $\delta<1-e^{-2L_0}$ we deduce
\begin{eqnarray*}
q(L_0) &\le& 2 \int_0^{A_\delta}
\frac{(1+u^2)g(u)}{\bfi(A_\delta
-L_0)-\bfi(A_\delta+L_0)}\,du
\\
&&{} +2\sqrt{2\pi}\int_{A_\delta}^{\pli}(u-L_0)
\bigl(1+u^2\bigr) \frac{e^{(u-L_0)^2/2}g(u)}{1-\delta-e^{-2L_0}}\,du
\\
&\le& C(A_\delta,L_0) + \frac{2e^{-L_0^2/2}}{1-\delta-e^{-2L_0}}\int
_{A_\delta}^{\pli} u\bigl(1+u^2
\bigr)e^{-L_0 u}\,du <\infty.
\end{eqnarray*}
Conclude that $\sup_{l,k}E_{f_0}^n|B_{lk}(\ix)| = O(1/n)$. Since
$\sum_{l,k}a_l^{-1} (\log a_l)^{-2\delta} <\infty$ the result follows.
\end{pf}
For the following theorem note that $\gamma=\delta=0$ gives $\|\cdot
\|_{0,2,0}=\|\cdot\|_2$.
%
\begin{theorem}\label{thm-l2}
With the notation of Theorem \ref{thm-mins}, suppose the product prior
$\Pi$ and $f_0$ satisfy Condition \ref{tolk}.
Then for any real numbers $\gamma,\delta$,
\[
E^{n}_{f_0} \int\|f-f_0\|_{\gamma,2,\delta}^2
\,d\Pi\bigl(f \given\mathbb{X}^{(n)}\bigr) = O \biggl(\sum
_{l, k}a_l^{2\gamma} (\log
a_l)^{-2\delta} \bigl(\sigma_l^2
\wedge n^{-1} \bigr) \biggr).
\]
\end{theorem}
\begin{pf}
We only prove $\gamma=\delta=0$; the general case is the same. With
the notation used in the proof of Theorem \ref{thm-mins}, using
Fubini's theorem,
\begin{eqnarray*}
E_{f_0}^{n} \int\|f-f_0\|_{2}^2\,d
\Pi\bigl(f \given\mathbb{X}^{(n)}\bigr) &=& \sum
_{l,k} E_{f_0}^{n} \int(\tlk-
\tolk)^2 \,d\Pi(\tlk\given\ix) \\
&=& \sum_{l,k}
E_{f_0}^{n} B_{lk}(\ix).
\end{eqnarray*}
In the proof of Theorem \ref{thm-mins}, the following two bounds have
been obtained, with the notation $ \cJ_n:=\{l\in\mathcal L, \rn
\sigma_l \ge S_0 \}$,
\[
\sup_{l\in\cJ_n, k} E_{f_0}^{n}
B_{lk}(\ix) = O\bigl( n^{-1}\bigr),\qquad 
\sup_{l \notin \mathcal J_n, k} E_{f_0}^n \sigma_{l}^{-2}
B_{lk}(\mathbb X) = O(1).
\]
%
For any $l\in\cJ_n^c$, by definition of $\cJ_n$ it holds $\sil
^2<S_0^2n^{-1}$, thus $\sil^2\le(1\vee S_0^2)( \sigma_l^2 \wedge
n^{-1})$. Similarly, if $l\in\cJ_n$, we have $n^{-1}\le(1\vee
S_0^{-2})( \sigma_l^2 \wedge n^{-1})$. 
\end{pf}
%
\begin{corollary} \label{postrat}
Set $\sil=|l|^{-{1}/{2}-\ga}$ or $\sil=2^{-({1}/{2}+\ga
)l}$ depending on the chosen $S$-regular basis of type either \textup{(a)} or
\textup{(b)}. Suppose that the conditions of Theorem \ref{thm-l2} are
satisfied. Then
\[
E^{n}_{f_0} \int\|f-f_0\|_{2}^2
\,d\Pi\bigl(f \given\mathbb{X}^{(n)}\bigr) = O \bigl(n^{-{2\ga}/({2\ga
+1})}
\bigr).
\]
\end{corollary}
\begin{pf}
For both types of basis $\sum_{l}|\mathcal Z_l|(\sigma_l^2 \wedge
n^{-1} )=O(n^{-{2\ga}/({2\ga+1})})$.
\end{pf}
%
\begin{remark} \label{rates}
The previous\vspace*{1pt} choice of $\sil$ entails a regularity condition on $f_0$
through condition (P1), namely $\sup_{k} |\tolk| \le M\sil$.
If $\sil=2^{-({1}/{2}+\ga)l}$ this amounts to the standard H\"
olderian condition if one uses a CDV wavelet basis, or a periodised
wavelet basis---any $f_0$ in $C^{\ga}$ from (\ref{hold}) satisfies
(P1) for such bases. For other bases similar remarks apply.
\end{remark}

\begin{corollary} \label{barrat}
Denote by $\bar{f}_n:=\bar{f}_n(\mathbb{X}^{(n)}):=\int fd\Pi(f
\given\mathbb{X}^{(n)})$
the posterior mean associated to the posterior distribution. Under the
conditions of Theorem \ref{thm-l2},
\[
E^{n}_{f_0} \|\bar{f}_n-f_0
\|_{2}^2 = O \biggl(\sum_{l, k}
\bigl( \sigma_l^2 \wedge n^{-1} \bigr) \biggr).
\]
\end{corollary}
\begin{pf}
Apply the Cauchy--Schwarz inequality and Theorem \ref{thm-l2}.
\end{pf}

\subsection{Convergence of the finite-dimensional distributions}
\label{fidi}

Consider again the posterior distribution $\Pi_n \equiv\Pi(\cdot
\given\mathbb{X}^{(n)})$ on $L^2$ from the beginning of this section (not
necessarily arising from a product measure). Let $V$ be any of the
finite-dimensional projection subspaces of $L^2$ defined in
Section \ref{sob}, equipped with the $L^2$-norm, and recall that $\pi
_V$ denotes the orthogonal $L^2$-projection onto $V$. For $z \in
H(\delta)$, define the transformation
\[
T_z \equiv T_{z,V}\dvtx  f \mapsto\sqrt n
\pi_{V}(f-z)
\]
from $H(\delta)$ to $V$, and consider the image measure $\Pi_n \circ
T_z^{-1}$. The finite-dimensional space $V$ carries a natural Lebesgue
product measure on it.

\begin{condition} \label{folk}
Suppose that $\Pi\circ\pi_{V}^{-1}$ has a Lebesgue-density $d\Pi_V$
in a neighborhood of $\pi_V(f_0)$ that is continuous and positive at
$\pi_V(f_0)$. Suppose also that for every $\delta>0$ there exists a
fixed $L^2$-norm ball $C=C_{\delta}$ in $V$ such that, for $n$ large enough,
$E_{f_0}^n (\Pi_n \circ T_{f_0}^{-1})(C^c) < \delta$.
\end{condition}

This condition requires that the projected prior has a continuous
density at $\pi_V(f_0)$ and that the image of the posterior
distribution under the finite-dimensional projection onto $V$
concentrates on a $1/\sqrt n$-neighborhood of $\pi_V(f_0)$. Let $\|
\cdot\|_{\mathrm{TV}}$ denote the total variation norm on the space of finite
signed measures on $V$, and $N(0, I)$ a standard Gaussian measure on $V$.

\begin{theorem} \label{proj}
Consider data generated from equation (\ref{wni}) under a fixed
function $f_0$, denote by $P_{f_0}^{n}$ the distribution of $\mathbb{X}^{(n)}$.
Assume Condition \ref{folk}. Then we have, as $n \to\infty$,
\[
\bigl\|\Pi_n \circ T_{\mathbb{X}^{(n)}}^{-1}-N(0, I)
\bigr\|_{\mathrm{TV}} \to^{P_{f_0}^n} 0.
\]
\end{theorem}

The proof of Theorem \ref{proj} is similar to the parametric proof in
Chapter 10 in \cite{V98}, and is omitted. In the special case of
product priors relevant for most examples in the present paper, one can
also derive the result directly from Theorem 1 in \cite{C12}: by
independence of the Gaussian coordinate experiments ${\langle}\psi_{lk},
\mathbb{X}^{(n)}{\rangle}\equiv\tolk+ \frac{1}{\rn}\elk$, when
estimating one or
more generally any finite number of the $\te_{lk}$'s, there is no loss
of information with respect to the case where all other $\te_{lk}$'s
would be known. Since the model is LAN with zero remainder, condition
(N) in \cite{C12} is satisfied, and condition (C) in \cite{C12}
amounts to asking that the full posterior concentrate at some rate
$\veps_n\to0$ in the $L^2$-norm (which for product priors is implied
by Corollary \ref{postrat}).


\subsection{\texorpdfstring{A BvM-theorem in $H(\delta)$}{A BvM-theorem in H(delta)}}

Let $\Pi_n=\Pi(\cdot|\mathbb{X}^{(n)})$ be the posterior
distribution on $L^2$.
Under the following Condition \ref{lolk}, which depends on a positive
real $\delta'$ to be specified in the sequel, we will prove that a
weak Bernstein--von Mises phenomenon holds true in $H(\delta)$ for any
$\delta>1/2$. For the product priors considered above we will then
verify Condition \ref{lolk} below.

\begin{condition} \label{lolk}
Suppose for every $\veps>0$ there exists a constant $0<M \equiv
M(\varepsilon)<\infty$ independent of $n$ such that, for any $n \ge
1$, some $\delta'>1/2$,
%
\begin{equation}
\label{bdclt} E_{f_0}^n \Pi\biggl[ \biggl\{ f\dvtx
\|f-f_0\|^2_{H(\delta')} > \frac
{M}{ n} \biggr\}
\Big| \mathbb{X}^{(n)} \biggr] \le\veps.
\end{equation}
Assume moreover that the conclusion of Theorem \ref{proj} holds true
for every $V$ (i.e., the finite-dimensional distributions converge).
\end{condition}


On $H(\delta)$ and for $z \in H(\delta)$, define the measurable map
\[
\tau_z\dvtx  f \mapsto\sqrt n (f-z).
\]
Recalling the definitions from Section \ref{sob}, consider $\Pi_n
\circ\tau_{\mathbb{X}^{(n)}}^{-1}$, a Borel probability measure on
$H(\delta)$.
Let $\mathcal N$ be the Gaussian measure on $H(\delta)$ constructed in
Section~\ref{sob} above.

\begin{theorem} \label{fbvm}
Fix $\delta>\delta'>1/2$, and assume Condition \ref{lolk} for such
$\delta'$. If $\beta$ is the bounded Lipschitz metric for weak
convergence of probability measures on $H(\delta)$, then as $n \to
\infty$, $\beta(\Pi_n \circ\tau_{\mathbb{X}^{(n)}}^{-1},\mathcal
N) \to0$ in
${P_{f_0}^n}$-probability.
\end{theorem}
\begin{pf}
It is enough to show that for every $\varepsilon>0$ there exists
$N=N(\varepsilon)$ large enough such that for all $n \ge N$,
\[
P_{f_0}^n \bigl(\beta\bigl(\Pi_n \circ
\tau_{\mathbb
{X}^{(n)}}^{-1},\mathcal N\bigr)>4\varepsilon\bigr) < 4
\varepsilon.
\]
Fix $\veps>0$, and let $V_J$ be the finite-dimensional subspace of
$L^2$ spanned by $\{\psi_{lk}\dvtx  k \in\mathcal Z_l, l \in\mathcal L,
|l|\le J\}$, for any integer $J\ge1$. Writing $\tilde\Pi_n$ for $\Pi
_n \circ\tau_{\mathbb{X}^{(n)}}^{-1}$ we see from the triangle inequality
\[
\beta(\tilde\Pi_n, \mathcal N) \le\beta\bigl(\tilde
\Pi_n, \tilde\Pi_n \circ\pi_{V_J}^{-1}
\bigr) + \beta\bigl(\tilde\Pi_n \circ\pi_{V_J}^{-1},
\mathcal N \circ\pi_{V_J}^{-1}\bigr) + \beta\bigl(\mathcal N
\circ\pi_{V_J}^{-1}, \mathcal N\bigr).
\]
The middle term converges to zero in $P_{f_0}^n$-probability for every
$V_J$, by convergence of the finite-dimensional distributions
(Condition \ref{lolk} and since the total variation distance dominates
$\beta$). Next we handle the first term.
Set $Q=M = M(\veps^2/4)$, and consider the random subset $D$ of
$H(\delta')$ defined as
\[
D = \bigl\{ g\dvtx  \|g+\dob\|_{H(\delta')}^2 \le Q \bigr\}.
\]
Under $P_{f_0}^n$we have $\tilde\Pi_n(D)=\Pi_n(D_n)$, where
\[
D_n=\bigl\{f\dvtx  \|f-f_0\|_{H(\delta')}^2
\le Q/n\bigr\}
\]
is the complement of the
set appearing in (\ref{bdclt}). In particular, using Condition~\ref
{lolk} and Markov's inequality yields
$P_{f_0}^n( \tilde\Pi_n(D^c)>\veps/4)\le\veps^2/\veps=\veps$.

If $Y_n\sim\tilde\Pi_n$ (conditional on $\mathbb{X}^{(n)}$), then
$\pi
_{V_J}(Y_n) \sim\tilde\Pi_n \circ\pi_{V_J}^{-1}$. For $F$ any
bounded function on $H(\delta)$ of Lipschitz-norm less than one,
\begin{eqnarray*}
&& \biggl\llvert\int_{H(\delta)} F \,d\tilde\Pi_n - \int
_{H(\delta)} F \,d\bigl(\tilde\Pi_n \circ
\pi_{V_J}^{-1}\bigr) \biggr\rrvert\\
&&\qquad= \bigl\llvert
E_{\tilde\Pi
_n} \bigl[F(Y_n)-F\bigl(\pi_{V_J}(Y_n)
\bigr) \bigr] \bigr\rrvert
\\
&&\qquad \le E_{\tilde\Pi_n} \bigl[ \bigl\|Y_n-\pi_{V_J}(Y_n)
\bigr\|_{H(\delta)} 1_{D}(Y_n) \bigr] + 2{\tilde
\Pi_n}\bigl(D^c\bigr),
\end{eqnarray*}
where $E_{\tilde\Pi_n}$ denotes expectation under $\tilde\Pi_n$
(given $\mathbb{X}^{(n)}$). With $y_{lk}=\langle Y_n, \psi_{lk}
\rangle$,
\begin{eqnarray*}
&&
E_{\tilde\Pi_n} \bigl[ \bigl\|Y_n-\pi_{V_J}(Y_n)
\bigr\|_{H(\delta)}^2 1_{D}(Y_n) \bigr]\\
&&\qquad=E_{\tilde\Pi_n} \biggl[ \sum_{l > J}
a_l^{-1} (\log a_l)^{-2\delta
} \sum
_k |y_{lk}|^2
1_{D}(Y_n) \biggr]
\\
&&\qquad= E_{\tilde\Pi_n} \biggl[ \sum_{l > J}
a_l^{-1} (\log a_l)^{2\delta'-2\delta-2\delta'} \sum
_k |y_{lk}|^2
1_{D}(Y_n) \biggr]
\\
&&\qquad\le (\log a_J)^{2\delta'-2\delta} E_{\tilde\Pi_n} \bigl[
\|Y_n\| _{H(\delta')}^2 1_{D}(Y_n)
\bigr] \\
&&\qquad\le2 (\log a_J)^{2\delta'-2\delta} \bigl[ Q + \|\dob
\|_{H(\delta
')}^2 \bigr].
\end{eqnarray*}
From the definition of $\be$ one deduces
\[
\be\bigl( \tilde\Pi_n, \tilde\Pi_n \circ
\pi_{V_J}^{-1} \bigr) \le2{\tilde\Pi_n}
\bigl(D^c\bigr) + \sqrt{2}(\log a_J)^{\delta'-\delta}
\sqrt{Q+ \|\dob\|^2_{H(\delta')}}.
\]
Since $a_J \to\infty$ as $J \to\infty$ we conclude that
$P_{f_0}^n(\be( \tilde\Pi_n, \tilde\Pi_n \circ\pi_{V_J}^{-1} )
>\veps)<2 \veps$ for $J$ large enough, combining the previous
deviation bound for $\tilde\Pi_n(D^c)$ and that $\|\dob\|_{H(\delta
')}$ is bounded in probability; cf. after (\ref{dual}) above.
A similar (though simpler) argument leads to
$ P_{f_0}^n(\beta(\mathcal N \circ\pi_{V_J}^{-1}, \mathcal N) >\veps
) <\veps$, using again that any random variable with law $\mathcal N$
has square integrable Hilbert-norm on $H(\delta')$. This completes the proof.
\end{pf}

\subsection{The BvM theorem for product priors} \label{pipp}

Combining Theorems \ref{thm-mins}, \ref{proj} and~\ref{fbvm}
implies that for product priors the weak Bernstein--von Mises theorem
in the sense of Definition \ref{bvmp} holds. The following results can
be seen to be uniform (``honest'') in all $f_0$ that satisfy Condition
\ref{tolk} with fixed constant $M$.

\begin{theorem} \label{pfbvm}
Suppose the assumptions of Theorem \ref{thm-mins} are satisfied and
that $\vphi$ is continuous near $\{\te_{0,lk}\}$ for every $k\in
\mathcal Z_l, l \in\mathcal L$. Let $\delta>1/2$. Then for $\beta$
the bounded Lipschitz metric for weak convergence of probability
measures on $H(\delta)$ we have, as $n \to\infty$, $\beta(\Pi_n
\circ\tau_{\mathbb{X}^{(n)}}^{-1},\mathcal N) \to^{P^n_{f_0}} 0$.
\end{theorem}
\begin{pf}
We only need to verify Condition \ref{lolk} with some $1/2<\delta
'<\delta$ so that we can apply Theorem \ref{fbvm}. From Theorem \ref
{thm-mins} with any such $\delta'$ in place of $\delta$, we see that
%
\begin{equation}
\label{mkv} 
nE_{f_0}^n\int\|f-f_0
\|^2_{H(\delta')}\,d\Pi\bigl(f|\mathbb{X}^{(n)}\bigr) =
O(1),
\end{equation}
which verifies the first part of Condition \ref{lolk} for some $M$
large enough using Markov's inequality. The second part follows from
verifying Condition \ref{folk} to invoke Theorem \ref{proj}: let $V$
be arbitrary. If $V_J$ is defined as in the proof of Theorem \ref
{fbvm}, and if $J$ is the smallest integer such that $V \subset V_J$, then
\[
\bigl\|\pi_V(f-f_0)\bigr\|_2^2 \le\bigl\|
\pi_{V_J}(f-f_0)\bigr\|_2^2 \le
a_J \log(a_J)^{2\delta'} \|f-f_0
\|^2_{H(\delta')}
\]
so that the second part of Condition \ref{folk} follows from the
estimate (\ref{mkv}) and again Markov's inequality, for $C$ a fixed
norm ball in $V$ of squared diameter of order $a_J \log(a_J)^{2\delta
'}M^2$. The first\vspace*{1pt} part of Condition \ref{folk} follows from the fact
that $\Pi\circ T^{-1}_{f_0}$ is a product measure in $V$ with bounded
marginals $\vphi_{lk}$ constant in $k$, and from the continuity
assumption on $\vphi$.
\end{pf}

\begin{theorem} \label{mom}
Suppose the assumptions of Theorem \ref{thm-mins} are satisfied and
that $\vphi$ is continuous near $\{\te_{0,lk}\}$ for every $k\in
\mathcal Z_l, l \in\mathcal L$. Let $\delta>1/2$ be arbitrary, let
$Y_n$ be a random variable drawn from $\Pi_n \circ\tau_{\mathbb
{X}^{(n)}}^{-1}$
(conditional on $\mathbb{X}^{(n)}$), and let $\bar f_n$ be the
(Bochner-) mean of
the posterior distribution $\Pi(\cdot|\mathbb{X}^{(n)})$. Then
$E[Y_n\given\mathbb{X}^{(n)}
] = \sqrt n(\bar f_n - \mathbb{X}^{(n)}) \to^{P_{f_0}^n} 0$ in
$H(\delta)$ as $n
\to\infty$.
\end{theorem}
\begin{pf}
Note that
\begin{eqnarray*}
E\bigl[\|Y_n\|_{H(\delta)}^2\given
\mathbb{X}^{(n)}\bigr] &=& \int\|h\| _{H(\delta)}^2\,d\Pi
_n \circ\tau_{\mathbb{X}^{(n)}}^{-1}(h)
\\
&\le& 2n \int\|f-f_0\|^2_{H(\delta)} \,d\Pi\bigl(f|
\mathbb{X}^{(n)}\bigr) + 2\| \dob\| _{H(\delta)}^2 \\
&=&
O_{P^n_{f_0}}(1)
\end{eqnarray*}
by Theorem \ref{thm-mins} and since $\|\dob\|_{H(\delta)}<\infty$
almost surely, as after (\ref{dual}). Moreover $Y_n \to N$ weakly in
$H(\delta)$ in $P_{f_0}^n$-probability, where $N \sim\mathcal N$, by
Theorem \ref{pfbvm}. By a standard uniform integrability argument
[using that $\{Y_n\dvtx  n \in\mathbb N\}$ has $H(\delta)$-norms with
uniformly bounded second moments and converges to $N$ weakly], and
arguing as in the last paragraph of Section \ref{weak} below, we
conclude $E[Y_n\given\mathbb{X}^{(n)}] \to EN$ in $H(\delta)$ in
$P_{f_0}^n$-probability, which implies the result since $EN=0$.
\end{pf}

\section{Remaining proofs}

\mbox{}

\begin{pf*}{Proof of Theorem \ref{sct}}
By Corollary 6.8.5 in \cite{B98} the image measure \mbox{$\mathcal N \circ
(\|\cdot\|_{H(\delta)})^{-1}$} of $\mathcal N$ under the norm mapping
is absolutely continuous on $[0, \infty)$, so the mapping
\[
\Phi\dvtx t \mapsto\mathcal N\bigl(B(0,t)\bigr)= \mathcal N \circ\bigl(\|\cdot\|
_{H(\delta)}\bigr)^{-1}\bigl([0,t]\bigr)
\]
is uniformly continuous and increasing on $[0, \infty)$. In fact, the
mapping is strictly increasing on $[0, \infty)$: using the results on
pages 213--214 in \cite{VV08b}, it suffices to show that any shell $\{
f\dvtx  s<\|f\|_{H(\delta)}<t\}, s<t$, contains an element of the RKHS
$L^2$ of $\cN$, which is obvious as $L^2$ is dense in $H(\delta)$.
Thus $\Phi$ has a continuous inverse $\Phi^{-1}\dvtx  [0,1) \to[0,\infty
)$. Since $\Phi$ is uniformly continuous for every $\epsilon>0$,
there exists $\delta>0$ small enough such that $|\Phi(t+\delta)-\Phi
(t)|<\epsilon$ for every $t \in[0, \infty)$. Now
\[
\mathcal N\bigl(\partial_\delta B(0,t)\bigr) = \mathcal N \bigl(B(0,t+
\delta)\bigr) - \mathcal N\bigl(B(0, t-\delta)\bigr)= \bigl|\Phi(t+\delta
)-\Phi(t-
\delta)\bigr| < 2\epsilon%
\]
for $\delta>0$ small enough, independently of $t$. Using (\ref
{nunif}) below we deduce that the balls $\{B(0, t)\}_{0 \le t <\infty
}$ form a $\mathcal N$-uniformity class, and we can thus conclude from
Definition \ref{bvmp} and the results in Section \ref{weak} below that
\[
\sup_{0\le t < \infty}\bigl\llvert\Pi\bigl(f\dvtx  \bigl\|f-\mathbb{X}^{(n)}
\bigr\| _{H(\delta)} \le t/\sqrt n |\mathbb{X}^{(n)}\bigr) - \mathcal N
\bigl(B(0,t)\bigr)\bigr\rrvert\to0
\]
in $P_{f_0}^n$-probability, as $n \to\infty$. This combined with
(\ref{qut}) gives
\[
\mathcal N\bigl(B(0,R_n)\bigr) = \mathcal N\bigl(B(0,R_n)
\bigr) - \Pi\bigl(f\dvtx  \bigl\|f-\mathbb{X}^{(n)}\bigr\| _{H(\delta)} \le
R_n/\sqrt n |\mathbb{X}^{(n)}\bigr) + 1-\alpha,
\]
which converges to $1-\alpha$ as $n \to\infty$ in
$P_{f_0}^n$-probability, and thus, by the continuous mapping theorem,
%
\begin{equation}
\label{rnlim} R_n \to^{P_{f_0}^n} \Phi^{-1}(1-\alpha)
\end{equation}
as $n \to\infty$. Now using this last convergence in probability,
\begin{eqnarray*}
P^n_{f_0}(f_0 \in C_n) &=&
P^n_{f_0} \bigl(f_0 \in B\bigl(
\mathbb{X}^{(n)}, R_n/\sqrt n\bigr)\bigr) =
P^n_{f_0} \bigl(0 \in B(\dob, R_n)\bigr)
\\
&=& P^n_{f_0} \bigl(0 \in B\bigl(\dob,
\Phi^{-1}(1-\alpha)\bigr)\bigr) +o(1)
\\
&=& \mathcal N\bigl(B\bigl(0, \Phi^{-1}(1-\alpha)\bigr)\bigr) +o(1)
\\
&=& \Phi\bigl(\Phi^{-1}(1-\alpha)\bigr)+o(1) = 1-\alpha+o(1),
\end{eqnarray*}
which completes the proof of the first claim. The second claim follows
from the same arguments combined with $\|\bar f_n - \mathbb{X}^{(n)}\|
_{H} =
o_P(n^{-1/2})$ which implies
\[
P^n_{f_0} \bigl(f_0 \in B(\bar
f_n, R_n/\sqrt n)\bigr)- P^n_{f_0}
\bigl(f_0 \in B\bigl(\mathbb{X}^{(n)}, R_n/\sqrt
n\bigr)\bigr) \to^{P_{f_0}^n} 0
\]
as $n \to\infty$, regardless of whether $R_n$ is defined via the
centering $T_n=\mathbb{X}^{(n)}$ or $T_n=\bar f_n$; cf. (\ref{qut}).
\end{pf*}

\begin{pf*}{Proof of Corollary \ref{unifset}}
By Theorems \ref{pfbvm} and \ref{mom} this prior satisfies the weak
Bernstein--von Mises phenomenon in the sense of Definition \ref{bvmp},
as well as (\ref{asylin}). The proof of coverage of $C'_n$ is thus the
same as in Theorem \ref{sct} noting that by hypothesis on $f_0$ the
probability $P^n_{f_0}(f_0 \in C'_n)$ in question equals
\[
P^n_{f_0} \bigl(\|f_0\|_{\gamma, \infty} \le M,
\|f_0-\bar f_n \| _{H(\delta)} \le R_n/
\sqrt n \bigr) = P^n_{f_0} \bigl(f_0 \in B(\bar
f_n, R_n/\sqrt n)\bigr).
\]
To control $|C'_n|_2$, pick two arbitrary functions $f_1,f_2$ in
$C'_n$, and let $g=f_1-f_2$. Then by construction and (\ref{rnlim}),
\[
\|g\|_{\gamma, 2,1} \le c\|g\|_{\gamma, \infty} \le2cM,\qquad \|g\|
_{H(\delta)}=O_P\bigl(n^{-1/2}\bigr).
\]
Choosing $J_n$ such that $2^{J_n} \sim n^{1/(2\gamma+1)}$,
\begin{eqnarray*}
\|g\|_2^2 &=& \sum_{l \ge J_0-1}
\sum_{k=0}^{2^l-1} \bigl|\langle g, \psi
_{lk} \rangle\bigr|^2
\\
&=& \sum_{l=J_0-1}^{J_n-1} l^{2\delta}2^{l}
2^{-l} l^{-2\delta} \sum_k \bigl|\langle
g, \psi_{lk} \rangle\bigr|^2 + \sum
_{l=J_n}^{\infty} 2^{-2l\gamma}l^{2}
2^{2l\gamma}l^{-2} \sum_k\bigl|\langle g,
\psi_{lk} \rangle\bigr|^2
\\
&\leq& 2^{J_n} J_n^{2\delta} \|g\|^2_{H(\delta)}
+ 2^{-2J_n\gamma
}J_n^{2} \|g\|^2_{\gamma,2,1}
\\
&=& O_{P} \biggl(\frac{2^{J_n}J_n^{2\delta}}{n} + 2^{-2J_n\gamma
}J_n^2
\biggr) = O_{P}\bigl(n^{-2\gamma/(2\gamma+1)} (\log n) ^{\kappa}\bigr)
\end{eqnarray*}
with constants independent of $g$, implying the same bound for $|C_n'|^2_2$.
\end{pf*}

\begin{pf*}{Proof of Corollary \ref{gaussset}}
By Theorems \ref{pfbvm} and \ref{mom} this prior satisfies the weak
Bernstein--von Mises phenomenon in the sense of Definition \ref{bvmp},
as well as (\ref{asylin}). By (\ref{normest}) we have $\|f_0\|
_{\gamma, 2, 1} \le M_n + 2 \delta+o_P(1)$ and so
$P^n_{f_0} (f_0 \in C_n'') = P^n_{f_0} (f_0 \in B(\bar f_n, R_n/\sqrt
n)) + o(1)$.
The proof of asymptotic $1-\alpha$-coverage of $C''_n$ is thus the
same as in Theorem \ref{sct}. Likewise, (\ref{normest}) implies $\Pi
_n(C''_n) = 1-\alpha+ o_P(1)$. To control $|C''_n|_2$, pick two
arbitrary functions $f_1,f_2$ in $C''_n $ and let $g=f_1-f_2$. Then by
(\ref{normest}) we have $\|g\|_{\gamma, 2,1}=O(M_n)=O_P(1)$ and by
(\ref{rnlim}) also $\|g\|_{H(\delta)}= O_P(n^{-1/2})$. The rest of
the proof is the same as in the previous corollary.
\end{pf*}

\begin{pf*}{Proof of Theorem \ref{convt}}
Since $f_0 \in L^1 \cap H^s_2$, we see by Fourier inversion on the
circle, the Cauchy--Schwarz inequality, and our assumption on the
equivalent Sobolev norm that
\begin{eqnarray*}
\|f \ast f_0\|_\infty&\le&\sum_m
\bigl|\hat f(m)\bigr| \bigl(1+|m|\bigr)^{-s} \bigl(1+|m|\bigr)^{s} \bigl|\hat{f}_0(m)\bigr|
\\
&\le&\biggl(\sum_m \bigl|\hat f(m)\bigr|^2
\bigl(1+|m|\bigr)^{-2s} \biggr)^{1/2} \biggl(\sum
_m \bigl|\hat f_0(m)\bigr|^2
\bigl(1+|m|\bigr)^{2s} \biggr)^{1/2} \\
&\le& C' \|f\|
_{H(\delta)}
\end{eqnarray*}
for any $\delta>0$, in particular $f\ast f_0$, for $f\in H(\delta),
f_0\in H_2^s$, defines a continuous function on $[0,1)$ (by Fourier
inversion), and the mapping $\lambda\dvtx  f \mapsto2f \ast f_0$ is linear
and continuous from $H(\delta)$ to $C([0,1))$; this argument is
adapted from Theorem 1 in \cite{N09}. By Definition \ref{bvmp} and
the continuous mapping theorem we thus have
$\beta((\Pi_n \circ\tau_{\mathbb{X}^{(n)}}^{-1}) \circ\lambda
^{-1},\mathcal
N \circ\lambda^{-1}) \to^{P_{f_0}} 0$
as $n \to\infty$, where $\beta$ is the bounded Lipschitz metric for
weak convergence in $C([0,1))$. Moreover from Corollary 6.8.5 in \cite
{B98} we deduce as in the proof of Theorem \ref{sct} that norm balls
$\{f\dvtx  \|f\|_\infty\le t\}_{0\le t < \infty}$ are $\mathcal N \circ
\lambda^{-1}$ uniformity classes for weak convergence, and that the
mapping $\Phi_\lambda\dvtx  t \mapsto\mathcal N \circ\lambda^{-1}(f\dvtx  \|
f\|_\infty\le t)$ from $[0, \infty)$ to $[0,1)$ is continuous and
increasing. In fact, it is strictly increasing, using the results on
pages 213--214 in \cite{VV08b} combined with the fact that the RKHS of
$\dob\ast f_0$, equal to $L^2\ast f_0$, contains functions of
arbitrary supremum norm.
Denote by $\Phi_\lambda^{-1}$ the continuous inverse of $\Phi
_\lambda$. As in the previous proofs, as $n \to\infty$
\[
\mathcal N \circ\lambda^{-1}\bigl(f\dvtx  \|f\|_\infty\le
R_n\bigr) - \bigl(\Pi_n \circ\lambda^{-1}\bigr)
\circ\theta_{\mathbb{X}^{(n)}\ast f_0}^{-1} \bigl(f\dvtx  \|f\| _\infty\le
R_n\bigr) \to^{P_{f_0}^n} 0,
\]
where $\theta_{\mathbb{X}^{(n)}\ast f_0}\dvtx  g \mapsto\sqrt n (g -
\mathbb{X}^{(n)}\ast f_0)$
maps $C([0,1)) \to C([0,1))$.

Thus, using the hypotheses on $\bar f_n$ and the posterior contraction
rate, the decomposition $f \ast f - g \ast g = 2(f-g)\ast g + (f-g)
\ast(f-g)$ and the convolution inequality $\|h \ast h'\|_\infty\le\|
h\|_2 \|h'\|_2$, we see
\begin{eqnarray*}
1-\alpha &=& \Pi_n \circ\kappa^{-1}\bigl(g\dvtx  \|g-\bar
f_n \ast\bar f_n\| _\infty\le R_n/
\sqrt n\bigr)
\\
&=& \Pi_n\bigl(f\dvtx \|f \ast f- \bar f_n \ast\bar
f_n\|_\infty\le R_n/\sqrt n\bigr)
\\
&\le&\Pi_n\bigl( f\dvtx  2\bigl\|\bigl(f-\mathbb{X}^{(n)}\bigr) \ast
f_0\bigr\|_\infty\le R_n/\sqrt n + r_n
\bigr) + o_P(1)
\\
&\le&\Pi_n\bigl( f\dvtx  2\rn\bigl\|\bigl(f-\mathbb{X}^{(n)}\bigr)
\ast f_0\bigr\|_\infty\le R_n + \delta_n
\bigr) + o_P(1) 
\end{eqnarray*}
with $\delta_n=r_n\sqrt{n}=o(1)$ as $n\to\infty$ by assumption.
Using the weak convergence property established above,
\[
1-\alpha\le\Phi_\la( R_n + \delta_n ) +
o_P(1)\quad \mbox{and similarly}\quad 1-\alpha\ge\Phi_\la(
R_n - \delta_n ) + o_P(1).
\]
From this we conclude $R_n \to^{P_{f_0}^n} \Phi_{\la}^{-1}(1-\alpha
)$ as $n \to\infty$. Now as above,
\begin{eqnarray*}
P^n_{f_0}(f_0 \ast f_0 \in
C_n) &=& P_{f_0}^n \bigl(\|f_0 \ast
f_0 - \bar f_n \ast\bar f_n\|_\infty
\le R_n /\sqrt n\bigr)
\\
&=& P_{f_0}^n\bigl(2\bigl\|(\bar f_n
-f_0) \ast f_0\bigr\|_\infty\le R_n/\sqrt
n\bigr) +o(1)
\\
&=& P_{f_0}^n\bigl(2 \sqrt n \bigl\|\bigl(\mathbb{X}^{(n)}-f_0
\bigr) \ast f_0\bigr\|_\infty\le\Phi_{\la
}^{-1}(1-
\alpha)\bigr) + o(1)
\\
&=& P_{f_0}^n\bigl(2\|\dob\ast f_0
\|_\infty\le\Phi_{\la}^{-1}(1-\alpha)\bigr) + o(1)
\\
&=& \Phi_{\la}\bigl( \Phi_{\la}^{-1}(1-\alpha)
\bigr) + o(1) = 1-\alpha+ o(1)
\end{eqnarray*}
completing the proof.
\end{pf*}

\begin{pf*}{Proof of Theorem \ref{linbvm}}
The proof is similar to the previous ones, using the continuous mapping
theorem for $L\dvtx  H(\delta) \to\mathbb R$, hence ommitted.
\end{pf*}

\begin{pf*}{Proof of Theorem \ref{bvm-func}}
The following notation is used in the proof:
\[
\te_n^* =\Psi(f_0) + \biggl\langle\dot{
\Psi}_{f_0}, \frac\dob\rn\biggr\rangle\quad\mbox{and}\quad \Phi_*(\cdot) = N
\bigl(0,\| \dot{\Psi}_{f_0} \|_2^2\bigr) \bigl((-
\infty, \cdot]\bigr).
\]
By definition of the quantile $\mu_n$ it holds
\begin{eqnarray*}
\frac\al2 &=& \Pi_n\circ\Psi^{-1}\bigl( (-\infty,
\mu_n] \bigr) = \Pi_n\bigl(\Psi(f) \le\mu_n\bigr)
\\
&=& \Pi_n \bigl(\Psi(f)- \Psi(f_0) \le\mu_n-
\Psi(f_0) \bigr)
\\
& =& \Pi_n \bigl( \bigl\langle\dot{\Psi}_{f_0},
f-\mathbb{X}^{(n)}\bigr\rangle\le\mu_n -
\te_n^* - \bigl[ \Psi(f)-\Psi(f_0) - \langle\dot{
\Psi}_{f_0}, f-f_0 \rangle\bigr] \bigr).
\end{eqnarray*}
The assumed contraction of the posterior in a $L^2$-neighborhood of
$f_0$ at rate $r_n$ together with (\ref{qfuz}) and the fact that $\rn
r_n=o(1)$ imply the existence of $\delta_n\to0$ such that
\begin{eqnarray*}
\frac\al2 &\le& \Pi_n\bigl( \rn\bigl\langle\dot{\Psi}_{f_0},
f-\mathbb{X}^{(n)} \bigr\rangle\le\rn\bigl(\mu_n -
\te_n^*\bigr) + \delta_n \bigr) + o_P(1),
\\
\frac\al2 &\ge& \Pi_n\bigl( \rn\bigl\langle\dot{\Psi}_{f_0},
f-\mathbb{X}^{(n)} \bigr\rangle\le\rn\bigl(\mu_n -
\te_n^*\bigr) - \delta_n \bigr) + o_P(1).
\end{eqnarray*}
Using the continuous mapping theorem and Definition \ref{bvmp},
\[
\mathcal\beta_\RR\bigl( \Pi_n \circ\ta_{\mathbb{X}^{(n)}}^{-1}
\circ(D\Psi_{f_0})^{-1}, \cN\circ(D\Psi_{f_0})^{-1}
\bigr) \to^{P_{f_0}^n} 0
\]
as $n \to\infty$. Note that $ \cN\circ(D\Psi_{f_0})^{-1} $ has
distribution function
$\Phi_*$.
Since the sets $\{(-\infty,t], t\in\RR\}$ form a uniformity class
for weak convergence towards a normal distribution, we obtain
\[
\frac\al2 \le\Phi_*\bigl( \rn\bigl(\mu_n - \te_n^*\bigr)
+ \delta_n \bigr) + o_P(1),\qquad\frac\al2 \ge\Phi_*\bigl(
\rn\bigl(\mu_n - \te_n^*\bigr) - \delta_n
\bigr) + o_P(1).
\]
%
From this we deduce $\mu_n = \te_n^* + \frac{1}{\rn}\Phi
_*^{-1}(\frac{\al}{2}) + o_P(1/\sqrt n)$. The quantile $\nu_n$
expands similarly, with $\Phi_*^{-1}(\frac{\al}{2})$ replaced by
$\Phi_*^{-1}(1-\frac{\al}{2})$.
By definition of $\theta^*_n$,
\begin{eqnarray*}
&& P_{f_0}^n \bigl(\Psi(f_0) \in(\mu_n,
\nu_n]\bigr)
\\
&&\qquad= P_{f_0}^n \biggl(\biggl\langle\dot{
\Psi}_{f_0}, \frac{\dob}{\rn} \biggr\rangle\in\biggl[
\frac{\Phi_*^{-1}(\alpha/2)}{\sqrt n} + o_P \biggl(\frac{1}{\sqrt n}
\biggr),\\
&&\qquad\hspace*{100.6pt}\frac{\Phi_*^{-1}(1-\alpha/2)}{\sqrt
n}+o_P \biggl(\frac{1}{\sqrt n} \biggr) \biggr]
\biggr)
\\
&&\qquad= P_{f_0}^n \bigl(\langle\dot{\Psi}_{f_0},
\mathbb W \rangle\in\bigl[\Phi_*^{-1}(\alpha/2), \Phi_*^{-1}(1-
\alpha/2)\bigr] \bigr) + o(1) \\
&&\qquad= 1-\alpha+o(1),
\end{eqnarray*}
completing the proof.
\end{pf*}

\subsection{Some weak convergence facts} \label{weak}

Let $\mu, \nu$ be Borel probability measures on a separable metric
space $(S,d)$. We call a family $\mathcal U$ of measurable real-valued
functions defined on $S$ a $\mu$-uniformity class for weak convergence
if for any sequence $\mu_n$ of Borel probability measures on $S$ that
converges weakly to $\mu$, we also have
%
\begin{equation}
\label{cunif} \sup_{u \in\mathcal U}\biggl\llvert\int
_S u(s) (d\mu_n-d\mu) (s)\biggr\rrvert\to0
\end{equation}
as $n \to\infty$. Necessary and sufficient conditions for classes
$\mathcal U$ of functions or sets $\{1_A\dvtx  A \in\mathcal A\}$ to form
uniformity classes are given in Billingsley and Tops{\o}e \cite
{BT67}. For any subset $A$ of $S$, define $A^\delta= \{x \in S\dvtx
d(x,A)<\delta\}$ and the $\delta$-boundary of $A$ by $\partial
_\delta A = \{x \in S\dvtx  d(x,A)<\delta, d(x,A^c)<\delta\}$. A family
$\mathcal A$ of measurable subsets of $S$ is a $\mu$-uniformity class
if and only if
%
\begin{equation}
\label{nunif} \lim_{\delta\to0} \sup_{A \in\mathcal A} \mu(
\partial_\delta A)=0;
\end{equation}
see Theorem 2 in \cite{BT67}. For classes of functions a similar
characterisation is available using moduli of continuity of the
involved functions; see Theorem 1 in \cite{BT67}. In particular the
bounded Lipschitz metric
\[
\beta(\mu, \nu) = \sup_{u \in \mathrm{BL}(1)}\biggl\llvert\int
_S u(s) (d\mu-d\nu) (s)\biggr\rrvert
\]
tests against the class
\[
\operatorname{BL}(1)= \Bigl\{f\dvtx  S \to\mathbb R, \sup_{s \in S}\bigl|f(s)\bigr| + \sup
_{s \ne
t, s,t \in S } \bigl|f(s)-f(t)\bigr|/d(s,t) \le1 \Bigr\},
\]
a uniformity class for any probability measure $\mu$. The metric
$\beta$ metrises weak convergence of probability measures on $S$
(\cite{D02}, Theorem 11.3.3).

We conclude with the following observation, which was used repeatedly
in our proofs: let $\mathcal P(S)$ denote the space of Borel
probability measures on $S$, let $(\Omega, \mathcal A, \mathbb P)$ be
a probability space, let $\mu_n\dvtx  (\Omega, \mathcal A, \mathbb P) \to
\mathcal P(S), n \in\mathbb N$, be random probability measures on $S$,
and let $\mu\in\mathcal P(S)$ be fixed. If $\beta(\mu_n, \mu) \to
^\mathbb P 0$ as $n \to\infty$, and if $\mathcal U$ is a $\mu
$-uniformity class, then the convergence in (\ref{cunif}) holds in
$\mathbb P$-probability, as is easily proved by contradiction and
passing to a.s. convergent subsequences. Likewise, if $(T,d')$ is a
metric space and $F\dvtx  S \to T$ a continuous mapping, then $\beta(\mu
_n \circ F^{-1}, \mu\circ F^{-1}) \to0$ in $\mathbb P$-probability.

\section*{Acknowledgement}

We would like to thank the Associate Editor for several fundamental
remarks that improved this article substantially.
I. Castillo would like to thank the Cambridge Statistical Laboratory as well
as Queens' College Cambridge for their hospitality during his visit to
Cambridge, where parts of this research were carried out.



\printaddresses

\end{document}